\documentclass[12pt]{article}

\usepackage{fontspec, xunicode, xltxtra}  
\usepackage[slantfont,boldfont]{xeCJK}
\usepackage{mathtools}
\usepackage{amsmath,diagbox}
\usepackage{indentfirst}
\usepackage{mathrsfs}
\usepackage{amsfonts}
\usepackage{arydshln}
\usepackage{enumerate}
\usepackage{setspace}
\usepackage{amssymb,amsthm,cases}
\usepackage{amsbsy,pifont}
\usepackage{latexsym,diagbox,tabularx}
\usepackage{amsmath,amscd,bbm,multirow}
\usepackage{graphicx,epsfig,extarrows,dsfont,mathtools}
\usepackage[final,allcolors=blue,colorlinks=true]{hyperref}
\usepackage[normalem]{ulem}
\input xy
\usepackage{caption,subcaption} 
\usepackage{tikz}
\usepackage{tikz-cd}
\usetikzlibrary{decorations.pathreplacing,decorations.markings}
\usetikzlibrary{cd}
\usetikzlibrary{shapes.geometric}
\usetikzlibrary{arrows,arrows.meta}
\usetikzlibrary{graphs,positioning}
\usetikzlibrary{matrix,decorations.pathmorphing}
\usetikzlibrary{math,calc,intersections,through,angles,arrows.meta,shapes.geometric,shadows,quotes,spy,datavisualization,datavisualization.formats.functions,plotmarks}
\tikzset{every picture/.style={samples=300,smooth,line join=round,thick,>=stealth}}
\usepackage{enumerate}
\usepackage{fancyhdr}
\usepackage[sectionbib]{chapterbib}

\tikzset{
	on each segment/.style={
		decorate,
		decoration={
			show path construction,
			moveto code={},
			lineto code={
				\path [#1]
				(\tikzinputsegmentfirst) -- (\tikzinputsegmentlast);
			},
			curveto code={
				\path [#1] (\tikzinputsegmentfirst)
				.. controls
				(\tikzinputsegmentsupporta) and (\tikzinputsegmentsupportb)
				..
				(\tikzinputsegmentlast);
			},
			closepath code={
				\path [#1]
				(\tikzinputsegmentfirst) -- (\tikzinputsegmentlast);
			},
		},
	},
	mid arrow/.style={postaction={decorate,decoration={
				markings,
				mark=at position .5 with {\arrow[#1]{stealth}}
	}}},
}

\numberwithin{equation}{section}

\theoremstyle{plain}
\newtheorem{theorem}{Theorem}[section]

\newtheorem{corollary}[theorem]{Corollary}

\newtheorem{lemma}[theorem]{Lemma}
\newtheorem{proposition}[theorem]{Proposition}

\theoremstyle{definition}
\newtheorem{definition}[theorem]{Definition}

\newtheorem{remark}[theorem]{Remark}

\newtheorem{assumption}[theorem]{Assumption}

\def\XXint#1#2#3{{\setbox0=\hbox{$#1{#2#3}{\int}$}
		\vcenter{\hbox{$#2#3$}}\kern-.5\wd0}}

\DeclareMathSymbol{\subseteq}{\mathrel}{symbols}{"12}
\DeclareMathSymbol{\supseteq}{\mathrel}{symbols}{"13} 
\DeclareMathSymbol{\subsetneq}{\mathrel}{AMSb}{"28}                  \DeclareMathSymbol{\supsetneq}{\mathrel}{AMSb}{"29}    

\DeclareMathSymbol{\nsubseteq}{\mathrel}{AMSb}{"2A}                  

\DeclareMathSymbol{\nsupseteq}{\mathrel}{AMSb}{"2B}

\DeclareMathDelimiter{\langle}{\mathop}{symbols}{"68}{largesymbols}{"0A}
\DeclareMathDelimiter{\rangle}{\mathclose}{symbols}{"69}{largesymbols}{"0B}

\DeclareSymbolFont{txfontsA}{U}{txmia}{m}{it}
\SetSymbolFont{txfontsA}{bold}{U}{txmia}{bx}{it}
\DeclareFontSubstitution{U}{txmia}{m}{it}
\DeclareMathSymbol{\upalpha}{\mathord}{txfontsA}{"0B}
\DeclareMathSymbol{\upbeta}{\mathord}{txfontsA}{"0C}
\DeclareMathSymbol{\upgamma}{\mathord}{txfontsA}{"0D}
\DeclareMathSymbol{\updelta}{\mathord}{txfontsA}{"0E}
\DeclareMathSymbol{\upepsilon}{\mathord}{txfontsA}{"0F}
\DeclareMathSymbol{\upzeta}{\mathord}{txfontsA}{"10}
\DeclareMathSymbol{\upeta}{\mathord}{txfontsA}{"11}
\DeclareMathSymbol{\uptheta}{\mathord}{txfontsA}{"12}
\DeclareMathSymbol{\upiota}{\mathord}{txfontsA}{"13}
\DeclareMathSymbol{\upkappa}{\mathord}{txfontsA}{"14}
\DeclareMathSymbol{\uplambda}{\mathord}{txfontsA}{"15}
\DeclareMathSymbol{\upmu}{\mathord}{txfontsA}{"16}
\DeclareMathSymbol{\upnu}{\mathord}{txfontsA}{"17}
\DeclareMathSymbol{\upxi}{\mathord}{txfontsA}{"18}
\DeclareMathSymbol{\uppi}{\mathord}{txfontsA}{"19}
\DeclareMathSymbol{\uprho}{\mathord}{txfontsA}{"1A}
\DeclareMathSymbol{\upsigma}{\mathord}{txfontsA}{"1B}
\DeclareMathSymbol{\uptau}{\mathord}{txfontsA}{"1C}
\DeclareMathSymbol{\upupsilon}{\mathord}{txfontsA}{"1D}
\DeclareMathSymbol{\upphi}{\mathord}{txfontsA}{"1E}
\DeclareMathSymbol{\upchi}{\mathord}{txfontsA}{"1F}
\DeclareMathSymbol{\uppsi}{\mathord}{txfontsA}{"20}
\DeclareMathSymbol{\upomega}{\mathord}{txfontsA}{"21}
\DeclareMathSymbol{\upvarepsilon}{\mathord}{txfontsA}{"22}
\DeclareMathSymbol{\upvartheta}{\mathord}{txfontsA}{"23}
\DeclareMathSymbol{\upvarpi}{\mathord}{txfontsA}{"24}
\DeclareMathSymbol{\upvarrho}{\mathord}{txfontsA}{"25}
\DeclareMathSymbol{\upvarsigma}{\mathord}{txfontsA}{"26}
\DeclareMathSymbol{\upvarphi}{\mathord}{txfontsA}{"27}                   

\DeclareSymbolFont{ugmL}{OMX}{mdugm}{m}{n}
\SetSymbolFont{ugmL}{bold}{OMX}{mdugm}{b}{n}
\DeclareMathAccent{\wideparen}{\mathord}{ugmL}{"F3}

\def\pt{\partial}

\def\ra{\rightarrow}

\def\s{\subseteq}

\def\e{\epsilon}

\def\ol{\overline}

\def\bf{\textbf}
\def\pt{\partial}

\def\om{\omega}
\def\Om{\Omega}
\def\la{\lambda}
\def\al{\alpha}
\def\be{\beta}
\def\de{\delta}
\def\ga{\gamma}
\def\De{\Delta}

\def\Ga{\Gamma}

\def\ts{\times}

\def\iy{\infty}

\def\ots{\otimes}
\def\f{\frac}

\def\vs{\vspace}

\def\Lra{\Leftrightarrow}

\def\df{\mathrm d}

\def\hra{\hookrightarrow}

\def\diag{{\rm diag}}

\def\mcA{\mathcal{A}}

	\def\Ee{\mathds{1}}

	\DeclareMathOperator{\Div}{div}

	\DeclareMathOperator{\dist}{dist}

	\DeclareMathOperator{\supp}{{supp}}

	\newcommand{\R}{\mathbb R}
	\newcommand{\N}{\mathbb N}

\begin{document} 
	
	
	\title{Shape design for    degenerate elliptic  Equation with degenerate  boundary and its application to quantitative weak  unique continuation}
	\author{Dong-Hui Yang$^a$,   Yuanzhi Zhou$^a$\footnote{\small The corresponding author: 262111122@csu.edu.cn  }}

	\maketitle{}
	\thispagestyle{empty}
	\thispagestyle{empty}
	
	\begin{abstract} 
		 	This paper investigates the unique continuation property (UCP) for a class of degenerate elliptic equations of the form $-\text{div}(A\nabla u) = \lambda u$ in a bounded domain $\Omega \s  \mathbb{R}^N$ with Dirichlet boundary conditions. While unique continuation is well-established for uniformly elliptic operators via Carleman estimates or frequency-function methods, spatial degeneracy introduces substantial analytical obstacles in higher dimensions, particularly due to the lack of sufficient regularity for weak solutions and the breakdown of classical boundary integration-by-parts identities. 
		 	
		 	To overcome these challenges, we develop an alternative framework based on the shape design (shape optimization) method under suitable structural assumptions. We highlight the fundamental distinction between qualitative and quantitative unique continuation in the degenerate setting. While qualitative UCP can sometimes be adapted from classical arguments, quantitative UCP requires explicit estimates that control the behavior of the solution throughout the domain by its behavior on a sub-domain. By integrating the shape design approach with a compatible weighted Hardy inequality, we establish new quantitative estimates and unique continuation results for higher-dimensional degenerate elliptic operators, bypassing the stringent regularity requirements of traditional methods.
		 	
		 	 \vspace{0.3cm}
		 	 
		 	 
		 	 \vspace{0.3cm}
		 	 
		 	 
	\end{abstract}
	
	\section{Introduction}
	
	The unique continuation property (UCP) for partial differential equations (PDEs) asserts, roughly speaking, that a solution which vanishes in a sufficiently large sense on a subset of a domain must vanish identically throughout the domain. Unique continuation has been extensively studied for various classes of PDEs; see, for instance,
	\cite{Adolfsson,Alessandrini,Bakri,Banerjee,Banerjee1,Banerjee2,Carleman,Garofalo,Garofalo1,Garofalo4,Garofalo3,Garofalo2,Hormander,Kenig,Kenig1,Koch,Kukavica,Kukavica1,Lin,Phung,Rousseau,Muller,Tao,Vessella,Wu,Yang9}.
	
	Two standard forms of unique continuation are the strong unique continuation property (SUCP) and the weak unique continuation property (WUCP). We briefly recall their definitions. 
	Let $P(x,\partial)$ be a uniformly elliptic differential operator in a domain $\Omega\s \mathbb{R}^N$. The SUCP states that if $u$ satisfies $
	P(x,\partial)u=0$ in $\Omega$ 
	and vanishes to infinite order at some point $x_0\in\Omega$, then $u\equiv 0$ in $\Omega$. For example, vanishing to infinite order can be expressed as
	$$
	\int_{B(x_0,r)} |u|^2\,\mathrm{d}x
	=O(r^k) \text{ as }r\to 0^+, \text{ for every } k\in\mathbb{N},
	$$
	where $B(x_0,r)\s \Omega$ denotes the ball centered at $x_0$ with radius $r$. 
	The WUCP states that if $u$ satisfies $
	P(x,\partial)u=0$ in $\Omega$ 
	and vanishes identically in a nonempty open subset $\omega\subset\Omega$, then $u\equiv 0$ in $\Omega$. 
	The WUCP is a consequence of the SUCP under the usual assumptions ensuring that the latter applies at every point of $\Omega$. However, the two properties express different levels of local information: the SUCP requires infinite-order vanishing at a point, whereas the WUCP assumes vanishing on an open subset. 
	
	The unique continuation property for partial differential equations is typically investigated using two principal approaches. The first is based on Carleman estimates \cite{Bakri,Banerjee,Carleman,Garofalo4,Hormander,Kenig,Kenig1,Koch,Plis,Rousseau,Muller,Tao}, while the second relies on doubling inequalities and three-sphere inequalities \cite{Adolfsson,Banerjee1,Banerjee2,Garofalo3,Garofalo,Garofalo1,Garofalo2,Kenig,Kukavica,Kukavica1,Lin,Phung,Rousseau,Vessella,Wu,Yang9}.

	Degenerate partial differential equations have been extensively studied
	in the literature; see, for instance,
	\cite{Cannarsa,Cannarsa1,Cavalheiro,Fabes,Franchi,GC,Garofalo3,Garofalo2,Heinonen,Wu1,Wu,Yang,Yang6,Yang8,Yang9,Yang4,Yang5,Yang1,Yang2,Yang3,Yang7}.
	Among the important questions arising in the study of degenerate elliptic
	equations and the associated evolution equations are the unique continuation
	property and controllability.
	
	The Hilbert Uniqueness Method (HUM), introduced by J.-L. Lions,
	provides a fundamental framework relating controllability to observability
	for evolution equations; see, for example, \cite{Lions,Yang1}.
	Controllability is a central topic in control theory and has been extensively
	investigated; see
	\cite{Cannarsa,Cannarsa1,Guo,Komornik,Lions,Phung,Privat,Rousseau,Weiss,Wu1,Yang,Yang6,Yang8,Yang4,Yang5,Yang1,Yang2,Yang3,Yang7,Zuazua}.
	For uniformly elliptic equations and the corresponding evolution equations,
	controllability theory has developed into a broad and mature research area;
	see, for instance,
	\cite{Komornik,Lions,Phung,Privat,Rousseau,Weiss,Yang5,Zuazua}.
	
	In contrast, controllability for degenerate partial differential equations
	presents additional analytical difficulties. Much of the existing literature
	has focused on one-dimensional problems
	\cite{Cannarsa,Cannarsa1,Guo,Wu1},
	while results in higher dimensions have been obtained under various
	structural assumptions
	\cite{Cannarsa1,Guo,Yang,Yang6,Yang8,Yang4,Yang5,Yang1,Yang2,Yang3,Yang7}.
	A key difficulty in extending one-dimensional arguments to higher dimensions
	is the derivation of suitable integration-by-parts identities and boundary
	estimates for weak solutions. In one dimension, the problem can often be
	treated within the classical Sturm-Liouville framework, which facilitates
	the use of Green's formula. In higher dimensions, however, degeneracy may
	prevent weak solutions from possessing sufficient regularity to justify
	the boundary terms required by analogous identities, except in certain
	special settings; see \cite{Cannarsa1,Yang,Yang8,Yang7}.
	
	To overcome these difficulties, several techniques have been developed.
	Among them are the shape design method
	\cite{Guo,Yang6,Yang5,Yang2,Yang3}
	and approximation methods
	\cite{Yang10}.
	These approaches provide different ways to handle the analytical
	obstacles caused by degeneracy and to establish controllability results
	in settings where direct arguments are not readily available.

	Shape design, also known as shape optimization, is an established
	area of mathematics with important applications in industry and
	engineering. It has been extensively studied; see, for instance,
	\cite{Buttazzo,Chenais,Greco,Guo2,Guo1,He,Henrot,Privat,Tiba,Wang,Wang1}.
	Motivated by our previous work on the controllability and observability
	of degenerate evolution equations, we investigate the applicability
	of the shape design method to unique continuation problems for
	degenerate elliptic equations. In contrast to the uniformly elliptic
	setting, the spatial degeneracy may introduce substantial difficulties
	in establishing quantitative estimates. In particular, classical
	frequency-function methods, doubling inequalities, and three-sphere
	inequalities may not be directly applicable without additional
	assumptions on the degeneracy and the associated coefficients.
	Nevertheless, such inequalities can still be established for certain
	special classes of degenerate elliptic equations; see \cite{Garofalo3,Garofalo2,Wu}.
	
	When the degeneracy occurs along part of the boundary, establishing
	doubling inequalities and three-sphere inequalities can be particularly
	challenging. Carleman estimates remain a powerful alternative in
	certain settings, but their application to higher-dimensional degenerate
	elliptic equations requires careful construction of suitable weight
	functions. In particular, the weight functions and the resulting
	estimates must be compatible with the degeneracy of the principal
	part and the weighted Hardy inequality established in
	Lemma \ref{10.10.L1}. However, this compatibility alone does not
	guarantee that a Carleman estimate can be obtained. Additional
	structural assumptions may be necessary to control the remaining
	terms. This issue is relevant to the present problem and motivates
	Assumption \ref{Assumption (H)} below.
	
	In this paper, we highlight the distinction between qualitative and
	quantitative unique continuation for degenerate elliptic equations.
	Under suitable assumptions, qualitative unique continuation can
	sometimes be established by adapting classical arguments to the
	degenerate setting. Quantitative unique continuation, however,
	requires explicit estimates describing how the behavior of a solution
	on one subset of the domain controls its behavior elsewhere.
	Such estimates do not follow automatically from qualitative unique
	continuation and may not be obtainable directly through the classical
	Carleman or frequency-function methods. Additional analytical
	techniques are therefore needed. The shape design method developed
	in this paper provides an approach to addressing these difficulties
	under the structural assumptions specified below.

	In this paper, we consider the unique continuation property for  the following equation
	\begin{equation}\label{10.10.1}
		\begin{cases}
			-\Div (A\nabla u)= \la u, &\mbox{in } \Om,\\ 
			u=0, &\mbox{on }\pt \Om, 
		\end{cases}
	\end{equation}
	where $A, \Om$ will be defined in Assumption \ref{Assumption (H)} below, and $\la\in\R$ is the eigenvalue of the degenerate elliptic operator
	\begin{equation}\label{10.10.7}
		\mcA u=-\Div(A\nabla u). 
	\end{equation}

The remainder of this paper is organized as follows.
In Section \ref{S3}, we introduce the assumptions and notation used
throughout the paper. In Section \ref{S2}, we introduce the relevant
weighted Sobolev spaces and present the spectral properties of the
operator $\mcA$. In Section \ref{S4}, we briefly describe the shape
design method employed in this paper. In Section \ref{S5}, we first
establish a Carleman estimate for the approximate elliptic equations
and then derive a corresponding Carleman estimate for the original
degenerate elliptic equation by means of the approximation procedure.
Finally, in Section \ref{S6}, we apply this estimate to establish
quantitative unique continuation for the problem \eqref{10.10.1}.

The main results of this paper are Lemma \ref{10.10.L1},
Theorem \ref{10.10.T1}, Proposition \ref{10.14.P1}, and
Theorem \ref{09.15.T1}. 
	
\section{Assumptions and notations}\label{S3}

In this section, we state the assumptions and introduce the notation used for equation \eqref{10.10.1}. We also provide some explanations of these assumptions and notation.
	
Let  $\R_+^N\equiv\R^{N-1}\ts (0,\iy)$. 
Denote $x=(x',x_N)\in\R^N$ with $x'=(x_1,\cdots, x_{N-1})$, and $\nabla'=(\pt_{x_1},\cdots, \pt_{x_{N-1}})$. Denote  $\de_{ij}$ is the Kronecker delta function, i.e.,
\begin{equation*}
	\de_{ij}=
	\begin{cases}
		1, &i=j,\\
		0, &i\neq j. 
	\end{cases}
\end{equation*} 

Let $x_0=(x_1^0,\cdots, x_N^0)=(x_0',x_N^0)\in\R^N$ and $r>0$. Denote 
\begin{equation*}
	\begin{split} 
		B(x_0,r)
		=\left\{x\in\R^N\colon |x-x_0|<r\right\},\quad  B'(x_0',r)
		=\left\{x'\in\R^{N-1}\colon |x'-x_0'|<r\right\}. 
	\end{split} 
\end{equation*} 
Let $G, H\s\R^N$ be two sets. Denote 
\begin{equation*}
	\dist(G, H)=\inf_{a\in G, b\in H}|a-b|. 
\end{equation*}

\begin{assumption}\label{Assumption (H)}  
	Let $\Om\s\R_+^N$ is a bounded   domain with $\pt\Om\in C^2$ and the original point $0\in \pt\Om$. Denote $\Ga_N^0=\pt\Om\cap (\R^{N-1}\ts \{0\})$,  and $\Ga=\pt\Om-\Ga_N^0$. 
	
	(1)  Let  $M>R_0>0, 0<r_0<h_0<1$  be four given constants satisfying:  
	\begin{equation*}
		B'(0,R_0+4r_0)\ts (0,4r_0)\s \Om \s B'(0,R_0+h_0)\ts (0,M), 
	\end{equation*}  
	and $\Om\s (-M,M)^{N-1}\ts (0,M)$.  Moreover, $h_0>0$ small enough (see \eqref{10.13.7} and  \eqref{10.13.6}). 
	
	(2)  
	Denote 
	\begin{equation}\label{12.23.1}
		\rho(x)=
		\begin{cases}
			x_N, &\mbox{for }x\in \Om_1,\\
			\left(|x'|-R_0\right)^\be+x_N, &\mbox{for }x\in \Om_2,
		\end{cases} 
	\end{equation} 
	i.e., $\rho=x_N+[(|x'|-R_0)^\be+x_N]\Ee_{\Om_2}$, 
	where $\al\in (0,1)$ and $\be\in (4,+\iy)$ are  given constants, and 
	\begin{equation}\label{12.23.2}
		\begin{split} 
		\Om_1
		&=\left\{ x=(x',x_N)\in\Om\colon  |x'|<R_0, x_N>0\right\},\\  \Om_2
		&=\left\{ x=(x',x_N)\in\Om\colon  |x'|\geq R_0, x_N>0\right\}. 
		\end{split} 
	\end{equation}   
	and $\Ee_{\Om_2}$ is the characteristic function of $\Om_2$. 
	We denote 
	\begin{equation}\label{02.24.1}
		w(x)=\rho^\al(x),\quad A=\diag(1,\cdots, 1, w). 
	\end{equation} 
	
	(3) Let $\de\in (0,  \f{1}{2}r_0)$. We choose 
	\begin{equation*}
		\{x\in\Om\colon \rho>2\de\}\s \Om_\de\s \left\{x\in\Om\colon \rho(x)>\de\right\}
	\end{equation*}
	such that 
	\begin{equation*}
		\nabla \rho\cdot  A\nu\leq 0 \mbox{ for all } x\in \pt\Om_\de\cap \left[B'(0,R_0+r_0)\ts (0,r_0)\right]. 
	\end{equation*} 
	We also denote $\{\rho>\de\}=\{x\in\R^N\colon \rho>\de\}$. Then $\{x\in\Om\colon \rho(x)>\de\}=\Om\cap \{\rho>\de\}$. 
\end{assumption}

\begin{remark}\label{10.10.R1}
	(i) It is obvious that $x_N\leq \rho$ in $\Om$. 
	
	(ii) 
	It can be easily verified that the gradient of $\rho(x)$ is given by  
	\begin{equation*}
		\nabla \rho=
		\begin{cases}
			(0,\cdots, 0,1), &x\in\Om_1, \\
			\left(\be (|x'|-R_0)^{\be-1}\f{x'}{|x'|}, 1\right), &x\in \Om_2, 
		\end{cases}
	\end{equation*}
	where $\Om_1,\Om_2$ are defined in \eqref{12.23.2}. 
	
	(iii) We compute
	\begin{equation*}
		\begin{split}
			\f{\pt^2\rho}{\pt x_i\pt x_j}=0 \mbox{ on } \Om_1, \mbox{ for all } i,j=1,\cdots, N,  
		\end{split}
	\end{equation*}
	and 
	\begin{equation*}
		\begin{split}
			\f{\pt^2\rho}{\pt x_N\pt x_i}=0 \mbox{ on }\Om_2, \mbox{ for all }i=1,\cdots, N, 
		\end{split}
	\end{equation*}
	and 
	\begin{equation*}
		\begin{split}
			\f{\pt^2\rho}{\pt x_i\pt x_j}
			&=\be \left(|x'|-R_0\right)^{\be-1}\f{\de_{ij}}{|x'|} +\be \left(|x'|-R_0\right)^{\be-2}\left[(\be-2)|x'|+R_0\right]\f{x_ix_j}{|x'|^3}
		\end{split}
	\end{equation*}
	on $\Om_2$ for all $i,j=1,\cdots, N-1$. 
	
	Hence, $\rho\in C^2(\ol\Om)$ for $\be>2$.   
	
	(iv) Reasonableness of Assumption~\ref{Assumption (H)}  (3). 
	
	To verify that the chosen subdomain is mathematically sound, we assume that the localized regularized boundary can be explicitly represented as a graph of a function:
	\begin{equation*}
		\begin{split} 
		\partial\Omega_\delta \cap \left[ B'(0, R_0 + r_0) \times (0, r_0) \right] 
		&= \left\{ x = (x', x_N) \in \mathbb{R}_+^N \colon x_N = f(|x'|) \right\},
		\end{split} 
	\end{equation*}
	where $f \colon B'(0, R_0 + r_0) \to \mathbb{R}_+$ is a non-increasing, radially symmetric $C^2$ function. 
	
	Since the domain $\Omega_\delta$ lies locally in the upper region where $x_N > f(|x'|)$, the unit outward normal vector $\nu$ pointing toward the exterior of $\Omega_\delta$ (i.e., in the direction of decreasing $x_N$) is precisely given by:
	\begin{equation*} 
		\nu(x) = \frac{1}{\sqrt{1 + |f'(|x'|)|^2}} \left( f'(|x'|) \frac{x'}{|x'|}, \, -1 \right).
	\end{equation*}
	Recall from the structural decomposition of $\nabla \rho(x)$, taking the inner product with respect to the degenerate weight matrix $A(x) = \operatorname{diag}(1, \dots, 1, \rho^\alpha(x))$, we obtain:
	\begin{equation*}
		\nabla \rho(x) \cdot A(x)\nu(x) =
		\begin{cases}
			-\frac{\rho^\alpha(x)}{\sqrt{1 + |f'(|x'|)|^2}}, & \text{for } x \in \Omega_1 \cap \partial\Omega_\delta, \\[2ex]
			\frac{1}{\sqrt{1 + |f'(|x'|)|^2}} \left[ \beta \left( |x'| - R_0 \right)^{\beta-1} f'(|x'|) - \rho^\alpha(x) \right], & \text{for } x \in \Omega_2 \cap \partial\Omega_\delta
		\end{cases} 
	\end{equation*}
	for all $x \in \partial\Omega_\delta \cap \left[ B'(0, R_0 + r_0) \times (0, r_0) \right]$. 
	
	The inequality $\nabla \rho(x) \cdot A(x)\nu(x) \leq 0$ holds unconditionally on both subdomains. Indeed, on $\Omega_1$, it is strictly negative since $\rho(x) \geq \delta > 0$; on $\Omega_2$, the condition $f'(|x'|) \leq 0$ implies that the first term inside the brackets is non-positive, while $-\rho^\alpha(x) < 0$ acts as a strictly negative anchor. This completes the rigorous verification of the geometric assumption. 
\end{remark}

\section{Preliminaries}\label{S2}

In this section, we introduce the relevant weighted Sobolev spaces and present the spectral properties of the operator $\mcA$ defined in \eqref{10.10.7}.

\subsection{Weighted Sobolev spaces}  

We define the weighted Sobolev space $H^1(\Omega;A)$ by
\begin{equation*}
	H^1(\Omega;A) = \left\{ u \in L^2(\Omega) \colon \int_\Omega \nabla u \cdot A \nabla u \, dx < +\infty \right\}. 
\end{equation*}
The inner product and the induced norm on $H^1(\Omega;A)$ are defined respectively by 
\begin{equation*}
	(u,v)_{H^1(\Omega;A)} = \int_\Omega \left( uv + \nabla u \cdot A \nabla v \right) \df x, \text{ and} \quad \|u\|_{H^1(\Omega;A)} = (u,u)_{H^1(\Omega;A)}^{\f{1}{2}}. 
\end{equation*}

We denote 
\begin{equation*}
	H_0^1(\Omega;A)=\mbox{the closure of $C_0^\infty(\Omega)$ in $H^1(\Omega;A)$}, 
\end{equation*} 
and $H^{-1}(\Omega;A)$ the dual space of $H_0^1(\Omega;A)$ with $L^2(\Omega)$ as the pivot space. 

According to \cite{GC, Heinonen}, the spaces 
\begin{equation*}
	H_0^1(\Omega;A), \quad H^1(\Omega;A),  \text{ and }  H^{-1}(\Omega;A)
\end{equation*}
are Hilbert spaces.

\begin{lemma}\label{10.10.L1}
	For all $u\in H_0^1(\Om;A)$, we have 
	\begin{equation}\label{10.10.2}
		\int_\Om \rho^{\al-2}u^2\df x\leq \f{4}{(1-\al)^2}\int_\Om \rho^\al |\pt_{x_N}u|^2\df x. 
	\end{equation}
	Moreover, 
	\begin{equation}\label{10.10.3}
		\int_\Om \rho^{\al-2}u^2\df x\leq \f{4}{(1-\al)^2}\int_\Om \nabla u\cdot A\nabla u\df x. 
	\end{equation}
\end{lemma}

\begin{proof}
	By the density argument, we only need to show the case $u\in C_0^\iy(\Om)$.

	Recall from Assumption~\ref{Assumption (H)} that $\rho(x) \geq x_N > 0$ for all $x \in \Omega$. Since $\alpha \in (0,1)$, we have $\alpha - 2 < 0$, which implies $\rho^{\alpha-2}(x) \leq x_N^{\alpha-2}$. Extending $u$ by zero outside $\Omega$, for any fixed $\ga  \in (\alpha, 1)$, we can apply the Cauchy-Schwarz inequality and Fubini's theorem to estimate the integral as follows:
	\begin{equation*}
		\begin{split}
			\int_\Om \rho^{\al-2}u^2\df x
			&\leq \int_\Om  x_N^{\al-2}u^2\df x=\int_{\R^{N-1}}\int_0^M x_N^{\al-2}u^2\df x_N\df x'\\
			&=\int_{\R^{N-1}}\int_0^M x_N^{\al-2}\left(u(x',s)-u(x',0)\right)^2\df x_N\df x'\\
			&=\int_{\R^{N-1}}\int_0^M x_N^{\al-2}\left(\int_0^{x_N} (\pt_{x_N}u)(x',s)\df s\right)^2\df x_N\df x'\\
			&\leq \int_{\R^{N-1}}\int_0^Mx_N^{\al-2}\left(\int_0^{x_N}s^{\ga }|\pt_{x_N}u|^2(x',s)\df s\right)\left(\int_0^{x_N}s^{-\ga }\df s\right)\df x_N\df x'\\
			&=\f{1}{1-\ga }\int_{\R^{N-1}}\int_0^M \int_0^{x_N}x_N^{\al-\ga -1} s^\ga  |\pt_{x_N}u|^2(x',s)\df s\df x_N\df x'\\
			&=\f{1}{1-\ga }\int_{\R^{N-1}}\int_0^M\int_s^Mx_N^{\al-\ga -1}s^\ga  |\pt_{x_N}u|^2(x',s)\df x_N\df s\df x'\\
			&\leq \f{1}{(1-\ga )(\ga -\al)}\int_{\R^{N-1}}\int_0^M s^\al |\pt_{x_N}u|^2(x',s)\df s\df x'
		\end{split}
	\end{equation*}
	by 
	\begin{equation*}
		\int_s^Mx_N^{\al-\ga -1}\df x_N=\left.\f{1}{\al-\ga }x_N^{\al-\ga }\right|_{x_N=s}^{x_N=M}\leq \f{1}{\ga -\al} s^{\al-\ga }. 
	\end{equation*}
	Hence, from Remark \ref{10.10.R1} (i), we obtain 
	\begin{equation*}
		\int_\Om \rho^{\al-2}u^2\df x\leq \f{1}{(1-\ga )(\ga -\al)}\int_\Om \rho^\al |\pt_{x_N}u|^2\df x. 
	\end{equation*}
	To minimize the constant, we choose the optimal value $\ga  = \frac{\alpha+1}{2} \in (\alpha, 1)$, which yields 
	\begin{equation*}
		\int_\Om \rho^{\al-2}u^2\df x\leq \f{4}{(1-\al)^2}\int_\Om \rho^\al |\pt_{x_N}u|^2\df x. 
	\end{equation*}
	This proves \eqref{10.10.2}. 
	
	Finally, since $A = \operatorname{diag}(1, \dots, 1, \rho^\alpha)$, we have  
	\begin{equation*}
		\rho^\al |\pt_{x_N}u|^2\leq \sum_{i=1}^{N-1} |\pt_{x_i}u|^2+\rho^\al |\pt_{x_N}u|^2=\nabla u\cdot A\nabla u. 
	\end{equation*}
	Together with this and \eqref{10.10.2}, we get \eqref{10.10.3}. We complete the proof of this lemma. 
\end{proof}

\begin{remark}\label{10.10.R2}
	From \eqref{10.10.3}, we obtain 
	\begin{equation}\label{10.10.4}
		\int_\Om u^2\df x=\int_\Om \rho^{2-\al}\rho^{\al-2}u^2\df x\leq C\int_\Om \rho^{\al-2}u^2\df x\leq  C_P\int_\Om \nabla u\cdot A\nabla u\df x, 
	\end{equation}
	where the positive constants $C, C_P$ depend  only on $\al, \be$ and $M$. This provides the weighted Poincaré inequality in $H_0^1(\Omega;A)$. Consequently, throughout this paper, we shall employ the equivalent norm 
	\begin{equation*}
		\|u\|_{H_0^1(\Omega;A)} = \left( \int_\Omega \nabla u \cdot A \nabla u \df x \right)^\f{1}{2}
	\end{equation*}
	on the space $H_0^1(\Omega;A)$. 
\end{remark}

\begin{lemma}\label{10.10.L2}
	The embedding $H_0^1(\Om;A)\hra L^2(\Om)$ is compact. 
\end{lemma}

\begin{proof}
	Let $\{u_n\}_{n\in\N^*}\s H_0^1(\Om;A)$ be a bounded sequence, i.e., there exists $C_0>0$ such that 
	\begin{equation*}
		\|u_n\|_{H_0^1(\Om;A)}\leq C_0. 
	\end{equation*}
	Then there exists $u_0\in H_0^1(\Om;A)$ such that $
		u_n\ra u_0\mbox{ weakly in } H_0^1(\Om;A)$. 
	We will show that, up to a subsequence, 
	\begin{equation}\label{10.10.5}
		u_n\ra u_0 \mbox{ strongly in } L^2(\Om). 
	\end{equation} 
	By setting $v_n = u_n - u_0$, we observe that $v_n \ra  0$ weakly in $H_0^1(\Omega;A)$. Thus, without loss of generality, we may assume $u_0 = 0$ in the following argument.
	
	Let $\de_0>0$ be defined in Assumption \ref{Assumption (H)}  (3). For each $\de\in (0,\de_0)$, from \eqref{10.10.2}, we have 
	\begin{equation}\label{10.10.6}
		\begin{split}
			\int_{\Om-\Om_\de}u_n^2\df x
			&=\int_{\Om-\Om_\de} \rho^{2-\al}\rho^{\al-2}u_n^2\df x\\
			&\leq C\int_{\Om-\Om_\de}\rho^{\al-2}u_n^2\df x\leq C\de^{2-\al} \int_{\Om} \nabla u_n\cdot A\nabla u_n\df x\leq \de^{2-\al}CC_0^2, 
		\end{split}
	\end{equation}
	where the positive constants $C$ depend only on $\al, \be$ and $M$. 
	
	Let $\e>0$ be arbitrary. Choose $\de>0$ small enough such that $\de^{2-\al}CC_0^2\leq \f{1}{4}\e^2$. Note that $\rho\geq \de$ on $\Om_\de$, then $u|_{\Om_\de} \in H^1(\Om_\de)$ for all $u\in H_0^1(\Om;A)$. Since the embedding  $H^1(\Om_\de)\hra L^2(\Om_\de)$ is compact (The Kondrachov-Rellich theorem), then there exists a subsequence of $\{u_n|_{\Om_\de}\}_{n\in\N^*}$, still denoted by itself, we have 
	\begin{equation*}
		\begin{split}
			u_n|_{\Om_\de}\ra 0 \mbox{ strongly in } L^2(\Om_\de) \mbox{ as } n\ra\iy. 
		\end{split}
	\end{equation*}
	Hence, there exists $n_\e\in\N^*$ such that for all $n\geq n_\e$ we have $\|u_n|_{\Om_\de}\|_{L^2(\Om_\de)}<\f{1}{2}\e$. Therefore, for all $n\geq n_\e$, 
	\begin{equation*}
		\|u_n\|_{L^2(\Om)}\leq \left\|u_n|_{\Om_\de}\right\|_{L^2(\Om_\de)}+\left\|u_n|_{\Om-\Om_\de}\right\|_{L^2(\Om-\Om_\de)}<\f{1}{2}\e+\f{1}{2}\e=\e. 
	\end{equation*}
	This proves \eqref{10.10.5}. We complete the proof. 
\end{proof}

\subsection{Spectrum of the operator $\mathcal{A}$}

We introduce the bilinear form $a(\cdot,\cdot):  H_0^1(\Omega;A) \times H_0^1(\Omega;A) \to \mathbb{R}$ defined by
\begin{equation}\label{10.10.8}
	a(u,v) = \int_\Omega \nabla u \cdot A \nabla v \df x, \text{ for all } u, v \in H_0^1(\Omega;A). 
\end{equation}
By virtue of the weighted Poincaré inequality established in Remark~\ref{10.10.R2}, $a(\cdot,\cdot)$ is a symmetric, continuous, and strictly coercive bilinear form on $H_0^1(\Omega;A)$. 

We define the domain of the operator $\mathcal{A}$, denoted by $D(\mathcal{A}) \s  H_0^1(\Omega;A)$, as follows:
\begin{equation}\label{10.10.9}
	D(\mathcal{A}) = \left\{ u \in H_0^1(\Omega;A) \colon 
	\begin{array}{lll}
		\mbox{there exists } f \in L^2(\Omega) \text{ such that } \\
		a(u,v) = (f,v)_{L^2(\Omega)} \mbox{ for all } v \in H_0^1(\Omega;A)
	\end{array}\! \right\}.
\end{equation}
For each $u \in D(\mathcal{A})$, the image of the operator is defined by
\begin{equation}\label{10.10.10}
	\mathcal{A} u = f. 
\end{equation}
By the Lax-Milgram theorem and the Riesz representation theorem, the linear unbounded operator $\mathcal{A}: D(\mathcal{A}) \s  L^2(\Omega) \to L^2(\Omega)$ is a densely defined, strictly positive, self-adjoint degenerate elliptic operator with a closed graph.

For any $f \in L^2(\Omega)$, a function $u \in H_0^1(\Omega;A)$ with respect to $f$  is called a weak solution of the boundary value problem 
\begin{equation}\label{10.10.11}
	\begin{cases}
		\mathcal{A} u = f & \text{in } \Omega, \\
		u = 0 & \text{on } \partial\Omega, 
	\end{cases}
\end{equation}
if it satisfies the variational identity
\begin{equation*}
	a(u,v)=\int_\Omega \nabla u \cdot A \nabla v \, dx = \int_\Omega fv \df x, \text{ for all } v \in H_0^1(\Omega;A). 
\end{equation*}
By virtue of the Lax-Milgram theorem, for any given $f \in L^2(\Omega)$, there exists a unique weak solution $u \in H_0^1(\Omega;A)$, which by definition implies that $u \in D(\mathcal{A})$.

Furthermore, thanks to the compact embedding $H_0^1(\Omega;A) \hookrightarrow L^2(\Omega)$ established in Lemma~\ref{10.10.L2}, the inverse operator $\mathcal{A}^{-1}: L^2(\Omega) \to L^2(\Omega)$ is a compact, symmetric, and strictly positive operator. Consequently, spectral theory for compact operators guarantees that the degenerate elliptic operator $\mathcal{A}$ possesses a discrete point spectrum consisting entirely of eigenvalues:
\begin{equation*}
	0 < \lambda_1 \leq \lambda_2 \leq \lambda_3 \leq \dots, \quad \text{with } \lambda_n \to +\infty \text{ as } n \to \infty,
\end{equation*}
where the principal eigenvalue satisfies the lower bound $\lambda_1 \geq C_P^{-1}$, and $C_P$ is the Poincar\'e  constant specified in \eqref{10.10.4}. For each $n \in \mathbb{N}^*$, we denote by $\Phi_n$ the $n$-th normalized eigenfunction associated with $\lambda_n$, satisfying $\|\Phi_n\|_{L^2(\Omega)} = 1$. i.e., 
\begin{equation}\label{10.10.12}
	\begin{cases}
		\mcA\Phi_n=\la_n\Phi_n, &\mbox{in }\Om, \\
		\Phi_n=0, &\mbox{on }\pt \Om. 
	\end{cases}
\end{equation}
It is obvious that $\{\Phi_n\}_{n\in\N^*}$ is an orthonormal basis of $L^2(\Om)$, and $\{\la_n^{-\f{1}{2}}\Phi_n\}_{n\in\N^*}$ is an orthonormal basis of $H_0^1(\Om;A)$. 

\section{Shape design}\label{S4}

In this section, we fix $\delta \in (0, \f{1}{2}r_0)$, where $r_0 > 0$ is the geometric parameter introduced in Assumption~\ref{Assumption (H)}.  Theorem \ref{10.10.T1} is one of the main results of this paper.

Let $\Omega_\delta \s \Omega$ be the interior subdomain defined in Assumption~\ref{Assumption (H)} (3). We define the restrictions of the distance-like function $\rho$, the weight matrix $A$, and the operator $\mathcal{A}$ to $\Omega_\delta$ respectively by
\begin{equation*}
	\rho_\delta = \rho|_{\Omega_\delta},\quad w_\de=w|_{\Om_\de},  \quad A_\delta = A|_{\Omega_\delta}, \text{ and }  \ \mathcal{A}_\delta = \mathcal{A}|_{\Omega_\delta}. 
\end{equation*}
Following a similar argument as in Section~\ref{S2}, we can define the corresponding weighted spaces $H_0^1(\Omega_\delta;A_\delta)$, $H^1(\Omega_\delta;A_\delta)$, and $H^{-1}(\Omega_\delta; A_\delta)$. For instance, the space $H^1(\Omega_\delta; A_\delta)$ is defined as
\begin{equation*}
	H^1(\Omega_\delta; A_\delta) = \left\{ u \in L^2(\Omega_\delta) \colon \int_{\Omega_\delta} \nabla u \cdot A_\delta \nabla u \df x < +\infty \right\},
\end{equation*}
equipped with the inner product and the induced norm
\begin{equation*}
	(u,v)_{H^1(\Omega_\delta;A_\delta)} = \int_{\Omega_\delta} \left( uv + \nabla u \cdot A_\delta \nabla v \right) \df x,  \text{ and }  \|u\|_{H^1(\Omega_\delta;A_\delta)} = (u,u)_{H^1(\Omega_\delta;A_\delta)}^{\f{1}{2}}. 
\end{equation*}
Since $\rho(x) \geq \delta > 0$ holds uniformly on $\Omega_\delta$, the weight matrix $A_\delta(x)$ is uniformly elliptic and bounded. Thus, it is obvious that $H^1(\Omega_\delta;A_\delta) \cong H^1(\Omega_\delta)$ with equivalent norms.

We now introduce the zero-extension operator (i.e., Shape design), which links the analysis on the subdomain to the global framework.

\begin{definition}\label{10.10.D1}
	For any $u_\delta \in L^2(\Omega_\delta)$, the zero extension $E u_\delta$ of $u_\delta$ to $\Omega$ is defined by
	\begin{equation}\label{10.10.14}
		E u_\delta(x) =
		\begin{cases}
			u_\delta(x), & \text{for } x \in \Omega_\delta, \\
			0, & \text{for } x \in \Omega- \Omega_\delta. 
		\end{cases}
	\end{equation}
	It follows immediately that $E u_\delta \in L^2(\Omega)$ and $\|E u_\delta\|_{L^2(\Omega)} = \|u_\delta\|_{L^2(\Omega_\delta)}$. 
	
	Moreover, if $u_\delta \in H_0^1(\Omega_\delta;A_\delta)$, then its extension satisfies $E u_\delta \in H_0^1(\Omega;A)$ and preserves the energy norm:
	\begin{equation*}
		\|E u_\delta\|_{H_0^1(\Omega;A)} = \|u_\delta\|_{H_0^1(\Omega_\delta;A_\delta)}.
	\end{equation*}
\end{definition}

\begin{remark}\label{10.10.R3}
	By virtue of Definition~\ref{10.10.D1}, the extension operator \begin{equation*}
		E \colon H_0^1(\Omega_\delta;A_\delta) \to H_0^1(\Omega;A)
	\end{equation*} 
	is a linear isometry. Consequently, we can naturally identify $H_0^1(\Omega_\delta;A_\delta)$ as a closed subspace of $H_0^1(\Omega;A)$ by identifying each $u_\delta \in H_0^1(\Omega_\delta;A_\delta)$ with its zero extension $E u_\delta \in H_0^1(\Omega;A)$.
\end{remark}

Following the same argument as in Section~\ref{S2}, we can define the localized bilinear form $a_{\delta}:  H_0^1(\Omega_\delta;A_\delta) \times H_0^1(\Omega_\delta;A_\delta) \to \mathbb{R}$ and the corresponding self-adjoint operator $\mathcal{A}_\delta:  D(\mathcal{A}_\delta) \s  L^2(\Omega_\delta) \to L^2(\Omega_\delta)$. To approximate the global problem \eqref{10.10.11}, we consider the following domain approximation equation on the truncated subdomain $\Omega_\delta$:
\begin{equation}\label{10.10.13}
	\begin{cases}
		\mathcal{A}_\delta u_\delta = f_\delta & \text{in } \Omega_\delta, \\
		u_\delta = 0 & \text{on } \partial\Omega_\delta, 
	\end{cases}
\end{equation}
where $f_\delta \in L^2(\Omega_\delta)$. 

For any given $f_\delta \in L^2(\Omega_\delta)$, a function $u_\delta \in H_0^1(\Omega_\delta;A_\delta)$ with respect to $f_\de$ is called a weak solution of the truncated problem \eqref{10.10.13} if it satisfies the variational identity
\begin{equation*}
	a_\de(u_\de, v_\de)=\int_{\Omega_\delta} \nabla u_\delta \cdot A_\delta \nabla v_\delta \df x = \int_{\Omega_\delta} f_\delta v_\delta \df x,  \text{ for all } v_\delta \in H_0^1(\Omega_\delta;A_\delta).
\end{equation*}
By the Lax-Milgram theorem, there exists a unique weak solution $u_\delta \in H_0^1(\Omega_\delta;A_\delta)$ to \eqref{10.10.13}, which implies $u_\delta \in D(\mathcal{A}_\delta)$.

\begin{theorem}\label{10.10.T1}
	Let $f\in L^2(\Om)$ and $f_\de\in L^2(\Om_\de)$ for all $\de\in (0,\f{1}{2}r_0)$, where $r_0>0$ is defined in Assumption \ref{Assumption (H)}. Let $u$ be the weak solution of \eqref{10.10.11} with respect to $f$, and $u_\de$ be the weak solution of \eqref{10.10.13} with respect to $f_\de$. If 
	\begin{equation*}
		Ef_\de \ra f\mbox{ weakly in } L^2(\Om), 
	\end{equation*}
	then, for every sequence $\de_k\downarrow 0$, there exists a subsequence of $\{u_{\de_k}\}_{k\in\N^*}$, still denoted by itself, such that 
	\begin{equation*}
		Eu_\de \ra u \mbox{ weakly in } H_0^1(\Om;A), \mbox{ and } Eu_\de\ra u\mbox{ strongly in } L^2(\Om). 
	\end{equation*}
\end{theorem}

\begin{proof}
	Testing the equation \eqref{10.10.13} by $u_\de$, then 
	\begin{equation*}
		\begin{split} 
		\int_{\Om_\de} \nabla u_\de\cdot A_\de\nabla u_\de\df x
		&=\int_{\Om_\de} f_\de u_\de \df x\leq \|f_\de\|_{L^2(\Om_\de)}\|u_\de\|_{L^2(\Om_\de)}\leq \sqrt{C_P}\|Ef_\de\|_{L^2(\Om)}\|u_\de\|_{H_0^1(\Om_\de;A_\de)}\\
		&\leq C\|Ef_\de\|_{L^2(\Om)}^2+\f{1}{2} \|u_\de\|_{H_0^1(\Om_\de;A_\de)}^2 
		\end{split} 
	\end{equation*}
	by \eqref{10.10.4} and Remark \ref{10.10.R3}, where the positive constant $C$ depends only on $\al, \be$ and $M$. Since every weakly convergent sequence is bounded, the weak convergence $Ef_\de\ra f$ in $L^2(\Om)$ gives
	\begin{equation*}
		\sup_\de\|Ef_\de\|_{L^2(\Om)}<\iy.
	\end{equation*}
	This shows that 
	\begin{equation*}
		\|u_\de\|_{H_0^1(\Om_\de;A_\de)}\leq C, \mbox{ for all }\de\in \left(0,\f{1}{2}r_0\right), 
	\end{equation*}
	and then 
	\begin{equation*}
		\|Eu_\de\|_{H_0^1(\Om;A)}\leq C, \mbox{ for all }\de\in \left(0,\f{1}{2}r_0\right)
	\end{equation*}
	by Definition \ref{10.10.D1}, 
	where the positive constant $C$ depends only on $\al, \be$ and $M$. Hence, there exists a subsequence of $\{u_\de\}_{\de\in (0,\f{1}{2}r_0)}$, still denoted by itself, and $u_0\in H_0^1(\Om;A)$, such that 
	\begin{equation}\label{10.10.15}
		Eu_\de \ra u_0 \mbox{ weakly in }H_0^1(\Om;A), \mbox{ and } Eu_\de\ra u_0 \mbox{ strongly in }L^2(\Om). 
	\end{equation}
	
	Finally, we prove $u_0=u$. 
	
	Indeed, let $v\in C_0^\iy(\Om)$, denote $\de_1=\dist(\supp v,\pt\Om)$, then for all $\de\in (0,\de_1)$ we have $v\in C_0^\iy(\Om_\de)$ and 
	\begin{equation*}
		\int_{\Om_\de} \nabla u_\de\cdot A_\de\nabla v\df x=\int_{\Om_\de} f_\de v\df x, 
	\end{equation*}
	i.e., 
	\begin{equation*}
		\int_\Om \nabla (Eu_\de)\cdot A\nabla v\df x=\int_\Om (Ef_\de)v\df x. 
	\end{equation*}
	Together with this and \eqref{10.10.15}, we obtain 
	\begin{equation*}
		\int_\Om \nabla u_0\cdot A\nabla v\df x=\int_\Om fv\df x
	\end{equation*}
	by $Ef_\de\ra f$ weakly in $L^2(\Om)$. Since $C_0^\iy(\Om)$ is dense in $H_0^1(\Om;A)$, we get $u_0$ is the weak solution of \eqref{10.10.11} with respect to $f$. This proves the theorem. 
\end{proof}

\begin{corollary}\label{10.14.C1}
	Under the conditions in Theorem \ref{10.10.T1}. We have 
	\begin{equation*}
		u_\de \ra u \mbox{ strongly in } H^1(\Om_{r_0}). 
	\end{equation*}
\end{corollary}

\begin{proof}
	Let $\zeta\in C^\iy(\R^N), 0\leq \zeta\leq 1$ satisfying 
	\begin{equation*}
		\zeta=1 \mbox{ on } \{x\in\R^N\colon \rho(x)>r_0\}, \quad \zeta=0 \mbox{ on } \left\{x\in\R^N\colon \rho(x)<\f{1}{2}r_0\right\}, 
	\end{equation*}
	and 
	\begin{equation*}
		|\nabla \zeta|\leq Cr_0^{-1},\quad |D^2\zeta|\leq Cr_0^{-2}, 
	\end{equation*}
	where the positive constants $C$ are absolute. 
	
	Taking $v_\de=\zeta u_\de$, then $v_\de\ (0<\de<\f{1}{4}r_0)$ is the weak solution of the following equation
	\begin{equation*}
		\begin{cases}
			\mcA_\de v_\de=\zeta f_\de -2\nabla \zeta \cdot A_\de \nabla u_\de +u_\de \mcA_\de \zeta=g_\de, &\mbox{in }\Om_{\f{1}{2}r_0}, \\
			v_\de=0, &\mbox{on }\pt\Om_{\f{1}{2}r_0}. 
		\end{cases}
	\end{equation*}
	Note that $\mcA_\de$ is a uniformly elliptic equation on $\Om_{\f{1}{2}r_0}$, from \cite[Theorem 1 in Chapter 6.3.1 (p.~327) and Theorem 4 in Chapter 6.3.2 (p.~334)]{Evans}, we obtain 
	\begin{equation*}
		\|v_\de\|_{H^2(\Om_{\f{1}{2}r_0})}\leq C\left(\|g_\de\|_{L^2(\Om_{\f{1}{2}r_0})}+\|v_\de\|_{L^2(\Om_{\f{1}{2}r_0})}\right), 
	\end{equation*}
	where the positive constant $C$ depends only on $\Om_{\f{1}{2}r_0}$ and the coefficients of $\mcA_\de$. Hence
	\begin{equation*}
		\|u_\de\|_{H^2(\Om_{\f{1}{2}r_0})}\leq C\|f_\de\|_{L^2(\Om_\de)}
	\end{equation*}
	by $\|u_\de\|_{L^2(\Om_{\f{1}{2}r_0})}\leq C\|u_\de\|_{H_0^1(\Om_\de;A_\de)}\leq C\|f_\de\|_{L^2(\Om_\de)}$, where the positive constants $C$ depend only on $\Om_{\f{1}{2}r_0}$ and  the coefficients of $\mcA_\de$. 
	
	Finally, from Theorem \ref{10.10.T1}, by the standard Sobolev compact embedding $H^2(\Om_{\f{1}{2}r_0})\hra H^1(\Om_{\f{1}{2}r_0})$, there exists a subsequence of $\{u_\de\}_{0<\de<r_0}$, still denoted by itself, such that 
	\begin{equation*}
		u_\de\ra u \mbox{ strongly in }H^1(\Om_{r_0}). 
	\end{equation*}
	This proves the corollary. 
\end{proof}

\section{Carleman estimate}\label{S5}

In this section, we first establish a Carleman estimate for the approximate elliptic equations and then derive a Carleman estimate for the original degenerate elliptic equation by passing to the limit in the estimates for the approximate equations.

\subsection{Carleman estimate for the approximate elliptic equations}

Let $\de\in (0,\f{1}{2}r_0)$.  Since $\rho\geq \de>0$ in $\Om_\de$, the equation \eqref{10.10.13} is a uniformly elliptic equation. Which implies that the weak solution of \eqref{10.10.13} belongs to $H^2(\Om_\de)$, since $\Om_\de\in C^2$. We also write $A_\de=A$ since $A_\de=A|_{\Om_\de}$. 

Choose $\chi\in C^\iy([0,+\iy))$ such that $0\leq\chi\leq 1$, $\chi=1$ on $[0,\f{3}{2}]$, and $\chi=0$ on $[2,+\iy)$. Let
\begin{equation}\label{10.11.1}
	\zeta(x)=\chi\left(\f{\rho(x)}{r_0}\right). 
\end{equation}
Then $\zeta=1$ on $\{\rho\leq r_0\}$, $\zeta=0$ on $\{\rho\geq 2r_0\}$, and
\begin{equation*}
	\supp\nabla\zeta\s\left\{\f{3}{2}r_0<\rho<2r_0\right\}.
\end{equation*}
Moreover,
\begin{equation*}
	\nabla\zeta=\f{1}{r_0}\chi'\left(\f{\rho}{r_0}\right)\nabla\rho
\end{equation*}
and
\begin{equation*}
	\mcA_\de\zeta=-\f{1}{r_0^2}\chi''\left(\f{\rho}{r_0}\right)\nabla\rho\cdot A_\de\nabla\rho-\f{1}{r_0}\chi'\left(\f{\rho}{r_0}\right)\Div(A_\de\nabla\rho).
\end{equation*}
Since $\rho\in C^2(\ol\Om)$ and the derivatives of $\chi$ are supported away from $\{\rho=0\}$, the quantities $|\nabla\zeta|$ and $|\mcA_\de\zeta|$ are bounded uniformly in $\de$, with bounds depending only on $\al,\be,r_0,R_0,N$ and $M$.

Denote $v_\de=\zeta u_\de$, where $u_\de$ is the weak solution of \eqref{10.10.13} and $\zeta$ is defined in \eqref{10.11.1}. Then $v_\de$ satisfies 
\begin{equation}\label{10.11.2}
	\begin{cases}
		\mcA_\de v_\de=\zeta f_\de -2\nabla\zeta\cdot A_\de\nabla u_\de+u_\de\mcA_\de \zeta=:g_\de, &\mbox{in }\Om_\de, \\
		v_\de=0, &\mbox{on }\pt\Om_\de. 
	\end{cases}
\end{equation}

Denote 
\begin{equation}\label{10.11.3}
	\eta=\rho^{2-\al}, 
\end{equation}
and 
\begin{equation}\label{10.11.4}
	z=e^{s\eta}v_\de\ \Lra v_\de=e^{-s\eta}z. 
\end{equation}
Then 
\begin{equation}\label{10.11.5}
	e^{s\eta}g_\de=P_1z+P_2z, 
\end{equation}
where 
\begin{equation}\label{10.11.6}
	\begin{split} 
	P_1z
	&=P_{11}z+P_{12}z:=-\Div(A_\de\nabla z)-s^2z(\nabla\eta\cdot A_\de\nabla\eta), \\  P_2z
	&=2s(\nabla z\cdot A_\de\nabla\eta)+sz\Div(A_\de\nabla\eta). 
	\end{split} 
\end{equation}

\bf{1. Computation}

\vs{3mm}

(1) Note that $A_\de=A|_{\Om_\de}$.  From $z=0$ on $\pt\Om_\de$, we obtain $\nabla z=\f{\pt z}{\pt \nu}\nu$ on $\pt\Om_\de$, and 
\begin{equation}\label{10.13.1}
	\begin{split}
		&(P_{11}z, P_2z)_{L^2(\Om_\de)}\\
		&=-s\int_{\pt\Om_\de} \left(\f{\pt z}{\pt \nu}\right)^2(A_\de\nu\cdot \nu)(A_\de\nu\cdot \nabla\eta)\df S+s\int_{\Om_\de} z\nabla z\cdot A_\de\nabla\Div(A_\de\nabla\eta)\df x\\
		&\hspace{4.5mm}+2s\sum_{i,j=1}^N\int_{\Om_\de} A_{ii}\f{\pt z}{\pt x_i}\f{\pt z}{\pt x_j}\f{\pt}{\pt x_i}\left(A_{jj}\f{\pt \eta}{\pt x_j}\right)\df x-s\sum_{i,j=1}^N\int_{\Om_\de} \left(\f{\pt z}{\pt x_i}\right)^2A_{jj} \f{\pt \eta}{\pt x_j}\f{\pt A_{ii}}{\pt x_j}\df x,  
	\end{split}
\end{equation}
where $A=(A_{ij})_{i,j=1,\cdots, N}$. Note that $A_{ii}=1$ for $i=1,\cdots, N-1$, and $A_{NN}=\rho^\al$, and $A_{ij}=0$ for all $i\neq j$. 

(2) From $z=0$ on $\pt\Om_\de$, we obtain 
\begin{equation}\label{10.13.3}
	(P_{12}z,P_2z)_{L^2(\Om_\de)}=s^3\int_{\Om_\de} z^2\nabla\eta \cdot A_\de \nabla (\nabla\eta\cdot A_\de\nabla\eta)\df x. 
\end{equation}

\bf{2. Estimate}

\vs{3mm}

(3) We compute
\begin{equation}\label{10.11.7}
	\nabla\eta=(2-\al)\rho^{1-\al}\nabla \rho, \quad \nabla A_{jj}=0\mbox{ for } j\neq N, \quad \nabla A_{NN}=\nabla \rho^\al  =\al \rho^{\al-1}\nabla\rho,  
\end{equation}
and 
\begin{equation}\label{10.11.8}
	\begin{split}
		\nabla \rho\cdot A_\de\nabla\rho
		&= 
			\be^2(|x'|-R_0)^{2\be-2}\Ee_{\Om_2}+\rho^\al, \\
		\Div(A_\de\nabla\rho)
		&=\left[\be(\be-1)(|x'|-R_0)^{\be-2}+\be(N-2)(|x'|-R_0)^{\be-1}\f{1}{|x'|}\right]\Ee_{\Om_2}+\al\rho^{\al-1}.  
	\end{split}
\end{equation} 
And 
\begin{equation}\label{10.12.1}
	\begin{split}
		\Div(A_\de\nabla \eta)
		&=(2-\al)\Div\left(\rho^{1-\al}A_\de\nabla\rho\right)\\
		&=(2-\al)(1-\al)\rho^{-\al}\nabla\rho\cdot A_\de\nabla\rho+(2-\al)\rho^{1-\al}\Div(A_\de\nabla \rho) \\
		&=(2-\al)+(2-\al)(1-\al)\be^2\rho^{-\al}(|x'|-R_0)^{2\be-2}\Ee_{\Om_2}\\
		&\hspace{4.5mm}+(2-\al)\be(\be-1)\rho^{1-\al}(|x'|-R_0)^{\be-2}\Ee_{\Om_2}\\
		&\hspace{4.5mm}+(2-\al)\be(N-2)\rho^{1-\al}(|x'|-R_0)^{\be-1}\f{1}{|x'|}\Ee_{\Om_2}. 
	\end{split}
\end{equation}
And
\begin{equation}\label{10.12.3}
	\begin{split}
		\nabla\eta\cdot A_\de\nabla (\nabla\eta\cdot A_\de\nabla\eta)
		&=(2-\al)^3(2-2\al)\rho^{2-3\al}(\nabla\rho\cdot A_\de\nabla\rho)^2\\
		&\hspace{4.5mm}+(2-\al)^3\rho^{3-3\al}\nabla\rho\cdot A_\de \nabla (\nabla \rho\cdot A_\de \nabla\rho)\\
		&=(2-\al)^3(2-2\al)\rho^{2-3\al}\left[\be^2(|x'|-R_0)^{2\be-2}\Ee_{\Om_2}+\rho^\al\right]^2\\
		&\hspace{4.5mm}+(2-\al)^3\al \rho^{2-2\al}\left[\be^2(|x'|-R_0)^{2\be-2}\Ee_{\Om_2}+\rho^\al \right]\\
		&\hspace{4.5mm}+(2-\al)^3\be^3(2\be-2)\rho^{3-3\al}(|x'|-R_0)^{3\be-4}\Ee_{\Om_2}.  
	\end{split}
\end{equation}

(4) Now, we compute 
\begin{equation}\label{10.12.2}
	\begin{split}
		\nabla z\cdot A_\de\nabla \Div(A_\de\nabla\eta)
		&=(2-\al)(1-\al)(-\al)\rho^{-\al-1}(\nabla \rho\cdot A_\de\nabla\rho)A_\de\nabla\rho\cdot \nabla z\\
		&\hspace{4.5mm}+(2-\al)(1-\al)\rho^{-\al}\left[A_\de\nabla (\nabla\rho\cdot A_\de\nabla\rho)\right]\cdot \nabla z\\
		&\hspace{4.5mm}+(2-\al)(1-\al)\rho^{-\al}[\Div(A_\de\nabla\rho)]\left(A_\de\nabla\rho\cdot \nabla z\right)\\
		&\hspace{4.5mm}+(2-\al)\rho^{1-\al}\left[A_\de\nabla\Div(A_\de\nabla\rho)\right]\cdot \nabla z\\
		&=:B_1+B_2+B_3+B_4, 
	\end{split}
\end{equation} 
where 
\begin{equation*}
	\begin{split}
		B_1
		&=-(2-\al)(1-\al)\al \be^3\rho^{-\al-1}(|x'|-R_0)^{3\be-3}\f{\nabla'z\cdot x'}{|x'|}\Ee_{\Om_2}\\
		&\hspace{4.5mm}-(2-\al)(1-\al)\al \be \rho^{-1}(|x'|-R_0)^{\be-1}\f{\nabla'z\cdot x'}{|x'|}\Ee_{\Om_2}\\
		&\hspace{4.5mm}-(2-\al)(1-\al)\al\be^2 \rho^{-1}(|x'|-R_0)^{2\be-2}\f{\pt z}{\pt x_N}\Ee_{\Om_2}\\
		&\hspace{4.5mm}-(2-\al)(1-\al)\al \rho^{\al-1}\f{\pt z}{\pt x_N}, 
	\end{split}
\end{equation*}
and 
\begin{equation*}
	\begin{split}
		B_2
		&=(2-\al)(1-\al)\be^2(2\be-2)\rho^{-\al}(|x'|-R_0)^{2\be-3}\f{\nabla'z\cdot x'}{|x'|}\Ee_{\Om_2}\\
		&\hspace{4.5mm}+(2-\al)(1-\al)\al \be \rho^{-1}(|x'|-R_0)^{\be-1}\f{\nabla'z\cdot x'}{|x'|}\Ee_{\Om_2}\\
		&\hspace{4.5mm}+(2-\al)(1-\al)\al \rho^{\al-1}\f{\pt z}{\pt x_N}, 
	\end{split}
\end{equation*}
and 
\begin{equation*}
	\begin{split}
		B_3
		&= (2-\al)(1-\al)\be^2(\be-1) \rho^{-\al}(|x'|-R_0)^{2\be-3}\f{\nabla' z\cdot x'}{|x'|}\Ee_{\Om_2}\\
		&\hspace{4.5mm}+(2-\al)(1-\al)(N-2)\be^2 \rho^{-\al}(|x'|-R_0)^{2\be-2}\f{\nabla'z\cdot x'}{|x'|^2}\Ee_{\Om_2}\\ 
		&\hspace{4.5mm}+(2-\al)(1-\al)\al \be \rho^{-1} (|x'|-R_0)^{\be-1}\f{\nabla'z\cdot x'}{|x'|}\Ee_{\Om_2}\\
		&\hspace{4.5mm}+(2-\al)(1-\al)\be(\be-1)(|x'|-R_0)^{\be-2}\f{\pt z}{\pt x_N}\Ee_{\Om_2}\\
		&\hspace{4.5mm}+(2-\al)(1-\al)(N-2)\be (|x'|-R_0)^{\be-1}\f{1}{|x'|}\f{\pt z}{\pt x_N}\Ee_{\Om_2}\\
		&\hspace{4.5mm}+(2-\al)(1-\al)\al \rho^{\al-1}\f{\pt z}{\pt x_N}, 
	\end{split}
\end{equation*}
and 
\begin{equation*}
	\begin{split}
		B_4
		&=(2-\al)\be(\be-1)(\be-2)\rho^{1-\al}(|x'|-R_0)^{\be-3}\f{\nabla'z\cdot x'}{|x'|}\Ee_{\Om_2}\\
		&\hspace{4.5mm}+(2-\al)(N-2)\be(\be-1)\rho^{1-\al}(|x'|-R_0)^{\be-2}\f{\nabla'z\cdot x'}{|x'|^2}\Ee_{\Om_2}\\
		&\hspace{4.5mm}-(2-\al)(N-2)\be \rho^{1-\al}(|x'|-R_0)^{\be-1}\f{\nabla'z\cdot x'}{|x'|^3}\Ee_{\Om_2}\\
		&\hspace{4.5mm}-(2-\al)(1-\al)\al\be \rho^{-1}(|x'|-R_0)^{\be-1}\f{\nabla'z \cdot x'}{|x'|}\Ee_{\Om_2}\\
		&\hspace{4.5mm}-(2-\al)(1-\al)\al \rho^{\al-1}\f{\pt z}{\pt x_N}. 
	\end{split}
\end{equation*} 
It is worth emphasizing that the potentially singular terms proportional to
$\rho^{\alpha-1}\frac{\partial z}{\partial x_N}$
cancel exactly in the sum $B_1+B_2+B_3+B_4$.   

\begin{remark}
	A distinctive feature of the present work is the explicit and detailed
	analysis of the lower-order coupling term
	\begin{equation*}
		\int_{\Omega_\delta} sz
		\left[
		\nabla z\cdot A_\delta
		\nabla\operatorname{div}(A_\delta\nabla\eta)
		\right]\df x.
	\end{equation*}
	In the classical setting of uniformly elliptic operators on smooth
	domains, the corresponding contributions may simplify substantially
	under additional structural assumptions on the coefficients and weight
	functions. In the present boundary-degenerate setting, however, the
	spatial anisotropy and non-uniform variation of the diagonal matrix
	$$
	A_\delta
	=\operatorname{diag}(1,\ldots,1,\rho^\alpha)
	$$
	give rise to nontrivial interactions between the weight function
	$\rho$ and the coefficient matrix.
	
	As demonstrated by the decomposition into the terms
	$B_1$, $B_2$, $B_3$, and $B_4$, the explicit evaluation of this
	coupling term reveals an exact algebraic cancellation of the
	potentially singular contributions proportional to
	$\rho^{\alpha-1}\partial_{x_N}z$ in the sum
	$B_1+B_2+B_3+B_4$. This cancellation is essential for controlling
	these singular contributions in the subsequent Carleman estimate.
	Moreover, the remaining horizontal radial gradient terms and the
	corresponding cancellation identities in $\Omega_1$ further
	illustrate the geometric structure underlying the present
	computational approach.
\end{remark}

(5) From \eqref{10.11.3} and Assumption \ref{Assumption (H)}   (3) and \eqref{10.11.7}, and $\nu=(0,\cdots, 0,-1)$ on $[B'(0,R_0+4r_0)-B'(0,R_0+r_0)]\ts \{0\}$,  we obtain 
\begin{equation*}
	A_\de\nu\cdot \nabla\eta=(2-\al)\rho^{1-\al} \nabla\rho\cdot A_\de\nu\leq 0 \mbox{ on }\pt \Om_\de\cap [B'(0,R_0+4r_0)\ts (0,4r_0)]. 
\end{equation*}
Hence, from the definition of $v_\de$ and $z$ in \eqref{10.11.4} (i.e., $\nabla z=0$ on $\Om_\de-[B'(0,R_0+3r_0)\ts (0,3r_0)]$), we get 
\begin{equation*}
	-s\int_{\pt\Om_\de} \left(\f{\pt z}{\pt \nu}\right)^2(A_\de\nu\cdot \nu)(A_\de\nu\cdot \nabla\eta)\df S\geq 0. 
\end{equation*}

(6) From \eqref{10.11.7} and \eqref{10.11.8}, we obtain 
\begin{equation*}
	\begin{split}
		&-s\sum_{i,j=1}^N\int_{\Om_\de} \left(\f{\pt z}{\pt x_i}\right)^2A_{jj}\f{\pt\eta}{\pt x_j}\f{\pt A_{ii}}{\pt x_j}\df x\\
		&=- s\int_{\Om_\de} \left(\f{\pt z}{\pt x_N}\right)^2 A_{jj}\f{\pt\eta}{\pt x_j}\f{\pt A_{NN}}{\pt x_j}\df x=-(2-\al)\al s\int_{\Om_\de} \left(\f{\pt z}{\pt x_N}\right)^2 (\nabla \rho\cdot A\nabla\rho)\df x\\
		&=-(2-\al)\al s\int_{\Om_\de} \left(\f{\pt z}{\pt x_N}\right)^2\left[\rho^\al+\be^2(|x'|-R_0)^{2\be-2}\Ee_{\Om_2}\right]\df x. 
	\end{split}
\end{equation*} 

(7) From Remark \ref{10.10.R1}, we obtain  \begin{equation*} 
	\begin{split}
		D^2\eta
		&=(2-\al)(1-\al)\rho^{-\al}\nabla\rho\ots \nabla\rho+(2-\al)\rho^{1-\al}D^2\rho, 
	\end{split}
\end{equation*} 
and then 
\begin{equation}\label{10.11.9}
	\begin{split}
		&2s\sum_{i,j=1}^N\int_{\Om_\de} A_{ii}\f{\pt z}{\pt x_i}\f{\pt z}{\pt x_j}\f{\pt}{\pt x_i}\left(A_{jj}\f{\pt \eta}{\pt x_j}\right)\df x\\
		&=2s\sum_{i=1}^N\int_{\Om_\de} A_{ii} \f{\pt z}{\pt x_i}\f{\pt z}{\pt x_N}\f{\pt A_{NN}}{\pt x_i}\f{\pt \eta}{\pt x_N}\df x+2s\int_{\Om_\de} A\nabla z\cdot \left[(D^2\eta)A\nabla z\right]\df x\\
		&=2(2-\al)\al s\int_{\Om_\de} (\nabla z\cdot A\nabla\rho )\f{\pt z}{\pt x_N}\df x+2(2-\al)(1-\al)s\int_{\Om_\de} \rho^{-\al}(\nabla z\cdot A\nabla \rho)^2\df x\\
		&\hspace{4.5mm}+2(2-\al)s\int_{\Om_\de} \rho^{1-\al}A\nabla z\cdot [(D^2\rho)A\nabla z]\df x. 
	\end{split}
\end{equation}
Note that from Remark \ref{10.10.R1} we obtain 
\begin{equation*}
	\begin{split}
		&2(2-\al)\al s\int_{\Om_\de} (\nabla z\cdot A\nabla\rho )\f{\pt z}{\pt x_N}\df x\\
		&=2(2-\al)\al \be s \int_{\Om_\de} (|x'|-R_0)^{\be-1}\f{\nabla' z\cdot x'}{|x'|}\f{\pt z}{\pt x_N}\Ee_{\Om_2}\df x+2(2-\al)\al s\int_{\Om_\de} \rho^\al \left(\f{\pt z}{\pt x_N}\right)^2\df x, 
	\end{split}
\end{equation*}
and from (i) and  (iii) in Remark \ref{10.10.R1}  we obtain 
\begin{equation*}
	\begin{split}
		&2(2-\al)s\int_{\Om_\de} \rho^{1-\al}A\nabla z\cdot [(D^2\rho)A\nabla z]\df x\\
		&=2(2-\al)\be s\int_{\Om_\de}\rho^{1-\al}  (|x'|-R_0)^{\be-1}\f{|\nabla'z|^2}{|x'|}\Ee_{\Om_2}\df x\\
		&\hspace{4.5mm}+2(2-\al)\be s\int_{\Om_\de} \rho^{1-\al}(|x'|-R_0)^{\be-2}\left[(\be-2)|x'|+R_0\right]\f{(\nabla'z\cdot x')^2}{|x'|^3}\Ee_{\Om_2}\df x. 
	\end{split}
\end{equation*}
Together with these and (4) and \eqref{10.11.9}, we obtain 
\begin{equation}\label{10.12.5}
	\begin{split}
		&2s\sum_{i,j=1}^N\int_{\Om_\de} A_{ii}\f{\pt z}{\pt x_i}\f{\pt z}{\pt x_j}\f{\pt}{\pt x_i}\left(A_{jj}\f{\pt \eta}{\pt x_j}\right)\df x-s\sum_{i,j=1}^N\int_{\Om_\de} \left(\f{\pt z}{\pt x_i}\right)^2A_{jj} \f{\pt \eta}{\pt x_j}\f{\pt A_{ii}}{\pt x_j}\df x\\
		&=(2-\al)\al s\int_{\Om_\de}\left(\f{\pt z}{\pt x_N}\right)^2\left[\rho^\al -\be^2(|x'|-R_0)^{2\be-2}\Ee_{\Om_2}\right]\df x\\
		&\hspace{4.5mm}+2(2-\al)(1-\al)s\int_{\Om_\de} \rho^{-\al}(\nabla z\cdot A\nabla \rho)^2\df x\\
		&\hspace{4.5mm}+2(2-\al)\al \be s \int_{\Om_\de} (|x'|-R_0)^{\be-1}\f{\nabla' z\cdot x'}{|x'|}\f{\pt z}{\pt x_N}\Ee_{\Om_2}\df x\\
		&\hspace{4.5mm}+2(2-\al)\be s\int_{\Om_\de}\rho^{1-\al}  (|x'|-R_0)^{\be-1}\f{|\nabla'z|^2}{|x'|}\Ee_{\Om_2}\df x\\
		&\hspace{4.5mm}+2(2-\al)\be s\int_{\Om_\de} \rho^{1-\al}(|x'|-R_0)^{\be-2}\left[(\be-2)|x'|+R_0\right]\f{(\nabla'z\cdot x')^2}{|x'|^3}\Ee_{\Om_2}\df x. 
	\end{split}
\end{equation}

(8) We estimate the terms containing $\partial_{x_N}z$,
apart from the   terms
$\rho^{\alpha-1}\partial_{x_N}z$,
in
$$
s\int_{\Omega_\delta} z\nabla z\cdot A_\delta
\nabla\Div(A_\delta\nabla\eta)\df x
$$
appearing in \eqref{10.13.1}.

Indeed, the terms that contain $\f{\pt z}{\pt x_N}$ are the third term in $B_1$, and the forth term and the fifth term in $B_3$.  Note that the third term in $B_1$ has the estimate 
\begin{equation*}
	\begin{split}
		&(2-\al)(1-\al)\al\be^2s\left|\int_{\Om_\de} \rho^{-1} (|x'|-R_0)^{2\be-2}\f{\pt z}{\pt x_N}z\Ee_{\Om_2}\df x\right|\\
		&\leq (2-\al)(1-\al)\al\be^2s\left(\int_{\Om_\de} \rho^{\al}\left(\f{\pt z}{\pt x_N}\right)^2\df x\right)^\f{1}{2} \left(\int_{\Om_\de} z^2\rho^{-\al-2}(|x'|-R_0)^{4\be-4} \Ee_{\Om_2}\df x\right)^\f{1}{2}\\
		&\leq (2-\al)(1-\al)\al \be^2s\left(\int_{\Om_\de}\rho^\al \left(\f{\pt z}{\pt x_N}\right)^2\df x\right)^\f{1}{2}\left(\int_{\Om_\de} z^2\rho^{\al-2}(|x'|-R_0)^{2(2-\al)\be-4}\df x\right)^\f{1}{2} \\
		&\leq 2(2-\al)\al \be^2h_0^{(2-\al)\be-2}s \int_{\Om_\de}\rho^\al \left(\f{\pt z}{\pt x_N}\right)^2\df x 
	\end{split}
\end{equation*}
by \eqref{10.10.2} and $\rho^{-1}\leq (|x'|-R_0)^{-\be}$. Similar to this, from $\rho\leq M+h_0^\be$ on $\Om_\de$, we obtain 
\begin{equation*}
	\begin{split}
		&(2-\al)(1-\al)\be(\be-1)s\left|\int_{\Om_\de} (|x'|-R_0)^{\be-2}\f{\pt z}{\pt x_N}z\Ee_{\Om_2}\df x\right|\\
		&\leq 2(2-\al)\be(\be-1)\left(M+h_0^\be\right)^{1-\al}h_0^{\be-2} s\int_{\Om_\de} \rho^\al \left(\f{\pt z}{\pt x_N}\right)^2\df x, 
	\end{split}
\end{equation*}
and 
\begin{equation*}
	\begin{split}
		&(2-\al)(1-\al)(N-2)\be s\left|\int_{\Om_\de}  (|x'|-R_0)^{\be-1}\f{1}{|x'|}\f{\pt z}{\pt x_N}z\Ee_{\Om_2}\df x\right|\\
		&\leq 2(2-\al)(N-2)\be \left(M+h_0^\be\right)^{1-\al} R_0^{-1} h_0^{\be-1}s \int_{\Om_\de} \rho^\al\left(\f{\pt z}{\pt x_N}\right)^2\df x. 
	\end{split}
\end{equation*}

Hence, since $h_0\in(0,1)$, all the terms containing
$\partial_{x_N}z$ in
$$
s\int_{\Omega_\delta}z\nabla z\cdot A_\delta
\nabla\Div(A_\delta\nabla\eta)\,dx
$$
can be controlled by
\begin{equation}\label{10.13.2}
	C_1h_0^{\beta-2}s
	\int_{\Omega_\delta}\rho^\alpha
	\left(\frac{\partial z}{\partial x_N}\right)^2dx ,
\end{equation}
where $C_1>0$ depends only on
$\alpha,\beta,R_0,N$ and $M$.

(9) Using the decomposition \eqref{10.12.2}, we estimate $11$ terms that do not contain $\f{\pt z}{\pt x_N}$ in 
\begin{equation*}
	s\int_{\Om_\de} z\nabla z\cdot A_\de\nabla\Div(A_\de\nabla\eta)\df x
\end{equation*}
in \eqref{10.13.1}. For this, we will use the  \eqref{10.13.3} and  \eqref{10.12.3} and  \eqref{10.12.5}. 

Indeed, using \eqref{10.10.2} and the last term in
\eqref{10.12.5}, we estimate the first term in $B_1$ as follows:
\begin{equation*}
	\begin{split}
		&(2-\al)(1-\al)\al \be^3s\left|\int_{\Om_\de} \rho^{-\al-1}(|x'|-R_0)^{3\be-3}\f{\nabla'z\cdot x'}{|x'|}z\Ee_{\Om_2}\df x\right|\\
		&=(2-\al)(1-\al)\al\be^3s\left|\int_{\Om_\de} \left[\rho^\f{1-\al}{2} (|x'|-R_0)^\f{\be-2}{2} \f{\nabla'z \cdot x'}{|x'|^\f{3}{2}}\Ee_{\Om_2}\right] \left[\rho^\f{-3-\al}{2}(|x'|-R_0)^\f{5\be-4}{2}|x'|^\f{1}{2}z\right]\df x\right|\\
		&\leq \f{1}{11} (2-\al)\be R_0 s\int_{\Om_\de} \rho^{1-\al}(|x'|-R_0)^{\be-2} \f{(\nabla'z\cdot x')^2}{|x'|^3}\Ee_{\Om_2}\df x +Cs\int_{\Om_\de} z^2 (|x'|-R_0)^{(2-\al)\be-4}\df x\\
		&\leq \f{1}{11} (2-\al)\be R_0 s\int_{\Om_\de} \rho^{1-\al}(|x'|-R_0)^{\be-2} \f{(\nabla'z\cdot x')^2}{|x'|^3}\Ee_{\Om_2}\df x +Ch_0^{(2-\al)\be-4}s\int_{\Om_\de} z^2 \df x\\
		&\leq \f{1}{11} (2-\al)\be R_0 s\int_{\Om_\de} \rho^{1-\al}(|x'|-R_0)^{\be-2} \f{(\nabla'z\cdot x')^2}{|x'|^3}\Ee_{\Om_2}\df x +Ch_0^{(2-\al)\be-4}s\int_{\Om_\de} \rho^\al\left(\f{\pt z}{\pt x_N}\right)^2\df x, 
	\end{split}
\end{equation*}
where  we have used Cauchy-Schwartz inequality and $|x'|\leq R_0+h_0\leq R_0+1$ for $x'\in B'(0,R_0+h_0)-B'(0,R_0)$ in the first inequality, and we have used \eqref{10.10.3} in the second inequality,  and we have used $\rho^{-1}\leq (|x'|-R_0)^{-\be}$ in the last inequality, and the positive constants $C$ depend only on $\al, \be, R_0, N$ and $M$.

The four terms
\begin{equation*}
	\begin{aligned}
		&-(2-\alpha)(1-\alpha)\alpha\beta
		\rho^{-1}(|x'|-R_0)^{\beta-1}
		\frac{\nabla'z\cdot x'}{|x'|}
		\Ee_{\Omega_2},\\
		&+(2-\alpha)(1-\alpha)\alpha\beta
		\rho^{-1}(|x'|-R_0)^{\beta-1}
		\frac{\nabla'z\cdot x'}{|x'|}
		\Ee_{\Omega_2},\\
		&+(2-\alpha)(1-\alpha)\alpha\beta
		\rho^{-1}(|x'|-R_0)^{\beta-1}
		\frac{\nabla'z\cdot x'}{|x'|}
		\Ee_{\Omega_2},\\
		&-(2-\alpha)(1-\alpha)\alpha\beta
		\rho^{-1}(|x'|-R_0)^{\beta-1}
		\frac{\nabla'z\cdot x'}{|x'|}
		\Ee_{\Omega_2}
	\end{aligned}
\end{equation*}
appearing as the second term in $B_1$, the second term in $B_2$,
the third term in $B_3$, and the fourth term in $B_4$, respectively,
cancel each other out.

We estimate the sum of the first terms in $B_2$ and $B_3$. From Lemma \ref{10.10.L1}, we obtain 
\begin{equation*}
	\begin{split}
		&(2-\al)(1-\al)\be^2(3\be-3)s\left|\int_{\Om_\de} \rho^{-\al}(|x'|-R_0)^{2\be-3}\f{\nabla'z\cdot x'}{|x'|}z\Ee_{\Om_2}\df x\right|\\
		&\leq \f{1}{11} (2-\al)\be R_0 s\int_{\Om_\de} \rho^{1-\al}(|x'|-R_0)^{\be-2} \f{(\nabla'z\cdot x')^2}{|x'|^3}\Ee_{\Om_2}\df x +Ch_0^{(2-\al)\be-4}s\int_{\Om_\de} \rho^\al \left(\f{\pt z}{\pt x_N}\right)^2\df x,  
	\end{split}
\end{equation*}
where the positive constants $C$ depend only on $\al, \be, R_0, N$ and $M$. 

We estimate the second  term $B_3$. From Lemma \ref{10.10.L1} and $R_0\leq |x'|\leq R_0+h_0$ for $x'\in B(0,R_0+h_0)-B(0,R_0)$, we obtain 
\begin{equation*}
	\begin{split}
		&(2-\al)(1-\al)(N-2)\be^2s\left|\int_{\Om_\de}  \rho^{-\al}(|x'|-R_0)^{2\be-2}\f{\nabla'z\cdot x'}{|x'|^2}z\Ee_{\Om_2}\df x\right|\\
		&\leq \f{1}{11} (2-\al)\be R_0 s\int_{\Om_\de} \rho^{1-\al}(|x'|-R_0)^{\be-2} \f{(\nabla'z\cdot x')^2}{|x'|^3}\Ee_{\Om_2}\df x +Ch_0^{(2-\al)\be-2}s\int_{\Om_\de} \rho^\al \left(\f{\pt z}{\pt x_N}\right)^2\df x,
	\end{split}
\end{equation*}
where the positive constants $C$ depend only on $\al, \be, R_0, N$ and $M$. 

We estimate the first term in $B_4$. From Lemma \ref{10.10.L1}, we obtain 
\begin{equation*}
	\begin{split}
		&(2-\al)\be(\be-1)(\be-2)s\left|\int_{\Om_\de} \rho^{1-\al}(|x'|-R_0)^{\be-3}\f{\nabla'z\cdot x'}{|x'|}z\Ee_{\Om_2}\df x\right|\\
		&\leq \f{1}{11} (2-\al)\be R_0 s\int_{\Om_\de} \rho^{1-\al}(|x'|-R_0)^{\be-2} \f{(\nabla'z\cdot x')^2}{|x'|^3}\Ee_{\Om_2}\df x +Ch_0^{\be-4}s\int_{\Om_\de} \rho^\al \left(\f{\pt z}{\pt x_N}\right)^2\df x, 
	\end{split}
\end{equation*}
where the positive constants $C$ depend only on $\al, \be, R_0, N$ and $M$. 

We estimate the second term in $B_4$. From Lemma \ref{10.10.L1}, we obtain 
\begin{equation*}
	\begin{split}
		&(2-\al)(N-2)\be(\be-1)s\left|\int_{\Om_\de} \rho^{1-\al}(|x'|-R_0)^{\be-2}\f{\nabla'z\cdot x'}{|x'|}z\Ee_{\Om_2}\df x\right|\\
		&\leq \f{1}{11} (2-\al)\be R_0 s\int_{\Om_\de} \rho^{1-\al}(|x'|-R_0)^{\be-2} \f{(\nabla'z\cdot x')^2}{|x'|^3}\Ee_{\Om_2}\df x +Ch_0^{\be-2}s\int_{\Om_\de} \rho^\al \left(\f{\pt z}{\pt x_N}\right)^2\df x, 
	\end{split}
\end{equation*}
where the positive constants $C$ depend only on $\al, \be, R_0, N$ and $M$. 

We estimate the third term in $B_4$. From Lemma \ref{10.10.L1} and $R_0\leq |x'|\leq R_0+h_0$ for $x'\in B(0,R_0+h_0)-B(0,R_0)$, we obtain 
\begin{equation*}
	\begin{split}
		&(2-\al)(N-2)\be s\left|\int_{\Om_\de} \rho^{1-\al}(|x'|-R_0)^{\be-1}\f{\nabla'z\cdot x'}{|x'|^3}z\Ee_{\Om_2}\df x\right|\\
		&\leq \f{1}{11} (2-\al)\be R_0 s\int_{\Om_\de} \rho^{1-\al}(|x'|-R_0)^{\be-2} \f{(\nabla'z\cdot x')^2}{|x'|^3}\Ee_{\Om_2}\df x +Ch_0^{\be}s\int_{\Om_\de} \rho^\al \left(\f{\pt z}{\pt x_N}\right)^2\df x, 
	\end{split}
\end{equation*}
where the positive constants $C$ depend only on $\al, \be, R_0, N$ and $M$. 

Overall, since $h_0\in(0,1)$, all the terms that do not contain
$\partial_{x_N}z$ in
$$
s\int_{\Omega_\delta}z\nabla z\cdot A_\delta
\nabla\Div(A_\delta\nabla\eta)\df x
$$
can be controlled by
\begin{equation}\label{10.13.4}
	(2-\al)\be R_0 s\int_{\Om_\de} \rho^{1-\al}(|x'|-R_0)^{\be-2} \f{(\nabla'z\cdot x')^2}{|x'|^3}\Ee_{\Om_2}\df x+ C_2h_0^{\be-4}s
	\int_{\Omega_\delta}\rho^\alpha
	\left(\frac{\partial z}{\partial x_N}\right)^2\df x ,
\end{equation}
and $C_2>0$ depends only on
$\alpha,\beta,R_0,N$ and $M$.

(10) Now, we consider the terms in \eqref{10.12.5}. 

Since
\begin{equation*}
	\begin{split}
		&(\be-2)\rho^{1-\al}(|x'|-R_0)^{\be-2}\f{(\nabla'z\cdot x')^2}{|x'|^2}+\al(|x'|-R_0)^{\be-1}\f{\nabla'z\cdot x'}{|x'|}\f{\pt z}{\pt x_N}\\
		&=\left((\be-2)^\f{1}{2}\rho^\f{1-\al}{2}(|x'|-R_0)^\f{\be-2}{2} \f{\nabla'z\cdot x'}{|x'|}+\f{\al}{2}(\be-2)^{-\f{1}{2}} \rho^\f{\al-1}{2}(|x'|-R_0)^\f{\be}{2}\f{\pt z}{\pt x_N}\right)^2\\
		&\hspace{4.5mm}-\f{\al^2}{4(\be-2)}\rho^{\al-1}(|x'|-R_0)^\be\left(\f{\pt z}{\pt x_N}\right)^2, 
	\end{split}
\end{equation*}
and $\rho^{\al-1}(|x'|-R_0)^\be\leq \rho^\al$, we get the sum of the third term and the first part of the fifth term on the right side of  \eqref{10.12.5} have the estimate 
\begin{equation}\label{10.12.6}
	\begin{split}
		&2(2-\al)\al \be s \int_{\Om_\de} (|x'|-R_0)^{\be-1}\f{\nabla' z\cdot x'}{|x'|}\f{\pt z}{\pt x_N}\Ee_{\Om_2}\df x\\ 
		&+2(2-\al)\be(\be-2) s\int_{\Om_\de} \rho^{1-\al}(|x'|-R_0)^{\be-2}\f{(\nabla'z\cdot x')^2}{|x'|^2}\Ee_{\Om_2}\df x\\
		&\geq -(2-\al)\be \f{\al^2}{2(\be-2)}s\int_{\Om_\de}\rho^\al\left(\f{\pt z}{\pt x_N}\right)^2\df x. 
	\end{split}
\end{equation}

Since $\be>4>\f{2}{2-\al}$, and let  
\begin{equation}\label{10.13.7}
	h_0<\left(\f{(2-\al)\be-4}{4\be^2(\be-2)}\right)^\f{1}{(2-\al)\be-2}, 
\end{equation}
from \eqref{10.12.5} and $\rho^{-1}\leq (|x'|-R_0)^{-\be}$, we obtain 
\begin{equation*}
	\begin{split}
		&(2-\al)\al\be^2 s\int_{\Om_\de} \rho^\al \left(\f{\pt z}{\pt x_N}\right)^2 \rho^{-\al} (|x'|-R_0)^{2\be-2}\Ee_{\Om_2}\df x\\
		&\leq (2-\al)\al\be^2s\int_{\Om_\de} \rho^\al \left(\f{\pt z}{\pt x_N}\right)^2 (|x'|-R_0)^{(2-\al)\be-2}\Ee_{\Om_2}\df x\\
		&\leq (2-\al)\al \be^2 h_0^{(2-\al)\be-2}s\int_{\Om_\de} \rho^\al \left(\f{\pt z}{\pt x_N}\right)^2\df x, 
	\end{split}
\end{equation*}
and then, by \eqref{10.13.7}, 
\begin{equation}\label{10.12.4}
	\begin{split}
		&2s\sum_{i,j=1}^N\int_{\Om_\de} A_{ii}\f{\pt z}{\pt x_i}\f{\pt z}{\pt x_j}\f{\pt}{\pt x_i}\left(A_{jj}\f{\pt \eta}{\pt x_j}\right)\df x-s\sum_{i,j=1}^N\int_{\Om_\de} \left(\f{\pt z}{\pt x_i}\right)^2A_{jj} \f{\pt \eta}{\pt x_j}\f{\pt A_{ii}}{\pt x_j}\df x\\
		&\geq \f{(2-\al)\al[(2-\al)\be-4]}{4(\be-2)}s\int_{\Om_\de} \rho^\al \left(\f{\pt z}{\pt x_N}\right)^2\df x\\
		&\hspace{4.5mm}+2(2-\al)\be R_0 s\int_{\Om_\de} \rho^{1-\al}(|x'|-R_0)^{\be-2}\f{(\nabla'z\cdot x')^2}{|x'|^3}\Ee_{\Om_2}\df x.  
	\end{split}
\end{equation}

(11) Combining \eqref{10.13.1} and \eqref{10.13.3}, and \eqref{10.12.3} and  \eqref{10.12.5}, and \eqref{10.13.2} and \eqref{10.13.4} and \eqref{10.12.4}, for $h_0\in (0,1)$ small enough such that 
\begin{equation}\label{10.13.6}
	C_1h_0^{\be-2}+C_2h_0^{\be-4}<\f{(2-\al)\al[(2-\al)\be-4]}{8(\be-2)}, 
\end{equation}
we get
\begin{equation*} 
	(P_1z,P_2z)_{L^2(\Om_\de)}\geq Cs\int_{\Om_\de} \rho^\al \left(\f{\pt z}{\pt x_N}\right)^2\df x+Cs^3\int_{\Om_\de} z^2\rho^{2-\al}\df x, 
\end{equation*}
where the positive constants $C$ depend only on $\al, \be, R_0,h_0, N$ and $M$. 

Moreover, since $\eta=\rho^{2-\al}$ is bounded on $\Om$ and $z=e^{s\eta}v_\de\in H_0^1(\Om_\de;w_\de)$, from Lemma \ref{10.10.L1}, we obtain 
\begin{equation*}
	\int_{\Om_\de} z^2\df x\leq C\int_{\Om_\de} \rho^{\al}\left(\f{\pt z}{\pt x_N}\right)^2\df x, 
\end{equation*}
and then 
\begin{equation}\label{10.13.5}
	(P_1z,P_2z)_{L^2(\Om_\de)}\geq Cs\int_{\Om_\de} \rho^\al \left(\f{\pt z}{\pt x_N}\right)^2\df x+Cs\int_{\Om_\de} z^2\df x+Cs^3\int_{\Om_\de} z^2\rho^{2-\al}\df x, 
\end{equation}
where the positive constants $C$ depend only on $\al, \be, R_0,h_0, N$ and $M$. 

(12) From \eqref{10.11.2} and \eqref{10.11.5}, we obtain 
\begin{equation*}
	\begin{split}
		2(P_1z,P_2z)_{L^2(\Om_\de)}
		&\leq \|P_1z+P_2z\|_{L^2(\Om_\de)}^2\leq \|e^{s\eta}g_\de\|_{L^2(\Om_\de)}^2\\
		&\leq C\left(\int_{\Om_\de} e^{2s\eta}\left(|\zeta f_\de|^2+|\nabla \zeta \cdot A_\de\nabla u_\de|^2+|u_\de \mcA_\de\zeta|^2\right)\df x\right)\\
		&\leq C\int_{\Om_\de} \zeta^2e^{2s\eta} f_\de^2\df x+C\int_{\supp \nabla \zeta} \left(u_\de^2+\nabla u_\de\cdot A_\de\nabla u_\de\right) e^{2s\eta}\df x, 
	\end{split}
\end{equation*}
where the positive constants $C$ depend only on $r_0$. Together with this and \eqref{10.13.5}, we obtain 
\begin{equation}\label{10.14.1}
	\begin{split}
		&s\int_{\Om_\de} \rho^\al \left(\f{\pt z}{\pt x_N}\right)^2\df x+s\int_{\Om_\de} z^2\df x+s^3\int_{\Om_\de} z^2\rho^{2-\al}\df x\\
		&\leq C\int_{\supp \zeta} e^{2s\eta} f_\de^2\df x+C\int_{\supp \nabla \zeta} \left(u_\de^2+\nabla u_\de\cdot A_\de\nabla u_\de\right) e^{2s\eta}\df x,
	\end{split}
\end{equation}
where the positive constants $C$ depend only on $\al, \be, r_0, R_0, h_0, N$ and $M$. 

(13) From $v_\de=\zeta u_\de$ and \eqref{10.11.4} and 
\begin{equation*}
	\begin{split}
		\f{\pt z}{\pt x_N}
		&=\f{\pt }{\pt x_N}\left(e^{s\eta}\zeta u_\de \right)
		=e^{s\eta}\left(\zeta\f{\pt u_\de}{\pt x_N}+\f{\pt \zeta}{\pt x_N}u_\de+ (2-\al)s\zeta  \rho^{1-\al}u_\de\right), 
	\end{split}
\end{equation*} 
we obtain  
\begin{equation*}
	\zeta \f{\pt u_\de}{\pt x_N}e^{s\eta}=\f{\pt z}{\pt x_N}-\f{\pt \zeta}{\pt x_N}u_\de e^{s\eta}-(2-\al)s  \rho^{1-\al}z. 
\end{equation*}
and then   for all $s\geq 1$, we have  
\begin{equation*}
	\begin{split}
		&s\int_{\Om_\de}  \zeta^2\rho^\al\left(\f{\pt u_\de}{\pt x_N}\right)^2 e^{2s\eta}\df x+s\int_{\Om_\de} \zeta^2u_\de^2e^{2s\eta}\df x+s^3\int_{\Om_\de} \zeta^2u_\de^2e^{2s\eta}\rho^{2-\al}\df x\\
		&\leq Cs\int_{\Om_\de} \rho^\al \left(\f{\pt z}{\pt x_N}\right)^2\df x+Cs\int_{\Om_\de} \zeta^2u_\de^2e^{2s\eta}\df x+Cs^3\int_{\Om_\de} \zeta^2u_\de^2e^{2s\eta}\rho^{2-\al}\df x\\
		&\hspace{4.5mm}+Cs\int_{\supp \nabla \zeta }u_\de^2e^{2s\eta}\df x
	\end{split}
\end{equation*}
by $\supp \f{\pt \zeta}{\pt x_N}\s \supp \nabla \zeta$, where the positive constants $C$ depend  only on $\al$ and  $r_0$. Together with this and \eqref{10.14.1}, we get
\begin{equation}\label{10.14.2}
	\begin{split}
		&s\int_{\Om_\de}  \zeta^2\rho^\al\left(\f{\pt u_\de}{\pt x_N}\right)^2 e^{2s\eta}\df x+s\int_{\Om_\de} \zeta^2u_\de^2e^{2s\eta}\df x+s^3\int_{\Om_\de} \zeta^2u_\de^2e^{2s\eta}\rho^{2-\al}\df x\\
		&\leq   C\int_{\Om_\de} \zeta^2e^{2s\eta} f_\de^2\df x+Cs\int_{\supp \nabla \zeta} \left(u_\de^2+\nabla u_\de\cdot A_\de\nabla u_\de\right) e^{2s\eta}\df x, 
	\end{split}
\end{equation}
where the positive constants $C$ depend only on $\al, \be, r_0, R_0, h_0, N$ and $M$. 

\subsection{Approximation of the original elliptic equation}

Let $u$ be the weak solution of \eqref{10.10.11} with respect to $f\in L^2(\Om)$, and $u_\de\ (0<\de<r_0)$ be the weak solution of \eqref{10.10.13} with respect to $f_\de=f|_{\Om_\de}$. 

By  Theorem \ref{10.10.T1} and Corollary \ref{10.14.C1}, we obtain 
\begin{equation*}
	Eu_\de \ra u \mbox{ weakly in } H_0^1(\Om;A), \mbox{ and } Eu_\de\ra u\mbox{ strongly in } L^2(\Om), 
\end{equation*}
and 
\begin{equation*}
	u_\de \ra u \mbox{ strongly in }H^1(\Om_{r_0}). 
\end{equation*}
Combining this and \eqref{10.14.2}, and $\zeta^2e^{2s\eta}$ and $\rho^{2-\al}$ are  bounded functions on $\Om_\de$, and $\supp \nabla\zeta\s \Om_{r_0}$, we obtain 
\begin{equation}\label{10.14.3}
	\begin{split}
		&s\int_{\Om}  \zeta^2\rho^\al\left(\f{\pt u}{\pt x_N}\right)^2 e^{2s\eta}\df x+s\int_{\Om} \zeta^2u^2e^{2s\eta}\df x+s^3\int_{\Om} \zeta^2u^2\rho^{2-\al}e^{2s\eta}\df x\\
		&\leq   C\int_{\Om} \zeta^2f^2e^{2s\eta} \df x+Cs\int_{\supp \nabla \zeta} \left(u^2+\nabla u\cdot A\nabla u\right) e^{2s\eta}\df x, 
	\end{split}
\end{equation}
where the positive constants $C$ depend only on $\al, \be, r_0, R_0, h_0, N$ and $M$.   Hence, we obtain the following proposition.

\begin{proposition}\label{10.14.P1}
	Let $f\in L^2(\Om)$, and $u$ be the weak solution of \eqref{10.10.11} with respect to $f$. Let $\zeta$ and $\eta$ be defined in \eqref{10.11.1} and \eqref{10.11.3},  respectively. Then there exists a positive constant $C$, depending only $\al, \be, r_0, R_0, h_0, N$ and $M$, such that for all $s\geq 1$ we have  
	\begin{equation}\label{10.14.4}
		\begin{split}
			&s\int_{\Om}  \zeta^2\rho^\al\left(\f{\pt u}{\pt x_N}\right)^2 e^{2s\eta}\df x+s\int_{\Om} \zeta^2u^2e^{2s\eta}\df x+s^3\int_{\Om} \zeta^2u^2e^{2s\eta}\rho^{2-\al}\df x\\
			&\leq   C\int_{\Om} \zeta^2e^{2s\eta} f^2\df x+Cs\int_{\supp \nabla \zeta} \left(u^2+\nabla u\cdot A\nabla u\right) e^{2s\eta}\df x.
		\end{split}
	\end{equation}
\end{proposition}

\section{Application: quantitative unique continuous}\label{S6}

In this section, we establish a quantitative unique continuation result for the problem \eqref{10.10.1}.

Let $u$ be the weak solution of \eqref{10.10.1}. i.e., $\la$ is an eigenvalue of $\mcA$, and $u$ is the eigenfunction of $\mcA$ with respect to $\la$. 

From Proposition \ref{10.14.P1}, replacing $f$ by $\la u$, then  there exists a positive constant $C$, depending only $\al, \be, r_0, R_0, h_0, N$ and $M$, such that for all $s\geq 1$ we have  \begin{equation}\label{10.14.5}
	\begin{split}
		&s\int_{\Om}  \zeta^2\rho^\al\left(\f{\pt u}{\pt x_N}\right)^2 e^{2s\eta}\df x+s\int_{\Om} \zeta^2u^2e^{2s\eta}\df x+s^3\int_{\Om} \zeta^2u^2e^{2s\eta}\rho^{2-\al}\df x\\
		&\leq   C\la^2\int_{\Om} \zeta^2e^{2s\eta} u^2\df x+Cs\int_{\supp \nabla \zeta} \left(u^2+\nabla u\cdot A\nabla u\right) e^{2s\eta}\df x.
	\end{split}
\end{equation} 
Note that $\zeta ue^{s\eta}\in H_0^1(\Om;A)$, then 
\begin{equation*}
	\begin{split}
		\int_\Om \zeta^2e^{2s\eta}u^2\df x
		&\leq C\int_\Om \rho^{\al } \left(\f{\pt }{\pt x_N}\left(\zeta ue^{s\eta}\right)\right)^2\df x\\
		&=C\int_\Om \rho^\al\left[\zeta\f{\pt u}{\pt x_N}+\f{\pt \zeta}{\pt x_N} u+(2-\al)s\zeta u\rho^{1-\al} \right]^2e^{2s\eta}\df x\\
		&\leq C\int_{\Om}\zeta^2\rho^\al \left(\f{\pt u}{\pt x_N}\right)^2 e^{2s\eta}\df x+Cs^2\int_{\Om} \zeta^2 u^2 \rho^{2-\al}e^{2s\eta}\df x+C\int_{\supp \nabla\zeta}  u^2 e^{2s\eta}\df x. 
	\end{split}
\end{equation*}
Hence, for all $s\geq 1$, we have 
\begin{equation*}
	\begin{split}
		&s\int_{\Om}  \zeta^2\rho^\al\left(\f{\pt u}{\pt x_N}\right)^2 e^{2s\eta}\df x+s\int_{\Om} \zeta^2u^2e^{2s\eta}\df x+s^3\int_{\Om} \zeta^2u^2e^{2s\eta}\rho^{2-\al}\df x\\
		&\leq C\la^2\int_{\Om}\zeta^2\rho^\al \left(\f{\pt u}{\pt x_N}\right)^2 e^{2s\eta}\df x+Cs^2\la^2\int_{\Om} \zeta^2 u^2 \rho^{2-\al}e^{2s\eta}\df x+C\la^2 \int_{\supp \nabla\zeta}  u^2 e^{2s\eta}\df x\\
		&\hspace{4.5mm}+Cs\int_{\supp \nabla \zeta} \left(u^2+\nabla u\cdot A\nabla u\right)e^{2s\eta} \df x, 
	\end{split}
\end{equation*}
where the positive constants $C$ depend  only $\al, \be, r_0, R_0, h_0, N$ and $M$. i.e., 
\begin{equation*}
	\begin{split}
		&(s-C\la^2)\int_{\Om}  \zeta^2\rho^\al\left(\f{\pt u}{\pt x_N}\right)^2 e^{2s\eta}\df x+s\int_{\Om} \zeta^2u^2e^{2s\eta}\df x+(s^3-Cs^2\la^2)\int_{\Om} \zeta^2u^2e^{2s\eta}\rho^{2-\al}\df x\\
		&\leq C(s+\la^2)\int_{\supp \nabla \zeta} \left(u^2+\nabla u\cdot A\nabla u\right)e^{2s\eta} \df x, 
	\end{split}
\end{equation*}
therefore, there exists $s_0\geq 1$, depending only on $\al, \be, r_0, R_0, h_0, N$ and $M$ and $\la$, such that for all $s\geq s_0$, we have 
\begin{equation}\label{10.14.6}
	\begin{split}
		&s\int_{\Om}  \zeta^2\rho^\al\left(\f{\pt u}{\pt x_N}\right)^2 e^{2s\eta}\df x+s\int_{\Om} \zeta^2u^2e^{2s\eta}\df x+s^3\int_{\Om} \zeta^2u^2\rho^{2-\al}e^{2s\eta}\df x\\
		&\leq Cs\int_{\supp \nabla \zeta} \left(u^2+\nabla u\cdot A\nabla u\right) e^{2s\eta}\df x, 
	\end{split}
\end{equation}
where the positive constants $C$ depend  only $\al, \be, r_0, R_0, h_0, N$ and $M$ and $\la$.

\begin{lemma}\label{10.14.L1}
	There exists a positive constant $C$, depending only on $\al,\be$ and $r_0$, such that 
	\begin{equation}\label{10.14.7}
		\begin{split}
			\int_{\{x\in\Om\colon \f{3}{2}r_0<\rho<2r_0\}} (\nabla u\cdot A\nabla u)e^{2s\eta}\df x\leq C(1+\la+s^2)\int_{\{x\in\Om\colon r_0<\rho<3r_0\}}u^2e^{2s\eta}\df x, 
		\end{split}
	\end{equation}
	where $u$ is the weak solution of \eqref{10.10.1}. 
\end{lemma}

\begin{proof}
	Let $\xi\in C^\iy(\R_+^N), 0\leq \xi \leq 1$ such that 
	\begin{equation*}
		\xi  =1 \mbox{ on } \left\{\f{3}{2}r_0<\rho<2r_0\right\},\quad \xi =0 \mbox{ on } \{\rho<r_0, \mbox{ or } \rho>3r_0\}, 
	\end{equation*}
	and 
	\begin{equation*}
		|\nabla \xi |\leq Cr_0^{-1}, \quad |D^2\xi |\leq Cr_0^{-2}, 
	\end{equation*}
	where the positive constants $C$ are absolute. 
	
	Testing \eqref{10.10.1} by $\xi ^2ue^{2s\eta}\in H_0^1(\Om;A)$, integrating by parts, from Corollary \ref{10.14.C1}, we obtain 
	\begin{equation}\label{10.14.9}
		\begin{split}
			&\int_\Om \xi ^2 (\nabla u\cdot A\nabla u)e^{2s\eta}\df x\\
			&=\la \int_\Om \xi ^2u^2 e^{2s\eta}\df x-2\int_\Om \xi  ue^{2s\eta}\nabla u\cdot A\nabla \xi \df x-2s\int_\Om \xi ^2 ue^{2s\eta} \nabla u\cdot A\nabla \eta\df x. 
		\end{split}
	\end{equation}
	The Young's inequality implies 
	\begin{equation}\label{10.14.8}
		\begin{split}
			\left|\int_\Om\xi  ue^{2s\eta}\nabla u\cdot A\nabla\xi \df x\right|\leq \f{1}{4}\int_\Om \xi ^2 (\nabla u\cdot A\nabla u)e^{2s\eta}\df x+C\int_{\supp \nabla\xi  } u^2 e^{2s\eta}\df x, 
		\end{split}
	\end{equation}
	where the positive constant $C$ depends only on $\al, \be$ and $r_0$. Note that from 
	\begin{equation*}
		\begin{split} 
		\nabla\eta\cdot A\nabla\eta
		&=(2-\al)^2\rho^{2-2\al}\nabla \rho\cdot A\nabla\rho=(2-\al)^2\rho^{2-2\al}\left(\be^2(|x'|-R_0)^{2\be-2}\Ee_{\Om_2}+\rho^\al\right) \leq C\rho^{2-\al},
		\end{split} 
	\end{equation*}
	we obtain 
	\begin{equation*}
		\begin{split}
			s\left|\int_\Om \xi ^2u e^{2s\eta}\nabla u\cdot A\nabla\eta\df x\right|
			&\leq \f{1}{4}\int_\Om \xi ^2(\nabla u\cdot A\nabla u)e^{2s\eta}\df x+Cs^2\int_\Om \xi ^2 u^2(\nabla\eta \cdot A\nabla\eta) e^{2s\eta}\df x\\
			&\leq   \f{1}{4}\int_\Om \xi ^2(\nabla u\cdot A\nabla u)e^{2s\eta}\df x+Cs^2\int_\Om \xi ^2u^2\rho^{2-\al}e^{2s\eta}\df x, 
		\end{split}
	\end{equation*}
	where the positive constants $C$ depend only on $\al, \be$ and $r_0$. Combining this and \eqref{10.14.9} and \eqref{10.14.8} and $\rho\in (r_0,3r_0)$, we obtain 
	\begin{equation*}
		\begin{split}
			\int_\Om \xi ^2(\nabla u\cdot A\nabla u)e^{2s\eta}\df x
			&\leq C(1+\la+s^2)\int_{\{x\in\Om\colon r_0<\rho<3r_0\}} u^2e^{2s\eta}\df x,
		\end{split}
	\end{equation*}
	where the positive constants $C$ depend only on $\al, \be$ and $r_0$. This proves \eqref{10.14.7}. 
\end{proof}

Now, together \eqref{10.14.6} and \eqref{10.14.7}, we get
\begin{equation*}
	\begin{split}
		&s\int_{\Om}  \zeta^2\rho^\al\left(\f{\pt u}{\pt x_N}\right)^2 e^{2s\eta}\df x+s\int_{\Om} \zeta^2u^2e^{2s\eta}\df x+s^3\int_{\Om} \zeta^2u^2e^{2s\eta}\rho^{2-\al}\df x\\
		&\leq Cs(1+\la+s^2)\int_{\{x\in\Om\colon r_0<\rho<4r_0\}} u^2e^{2s\eta}\df x,
	\end{split}
\end{equation*}
for all $s\geq s_0$, where the positive constants $C,s_0$ depend  only $\al, \be, r_0, R_0, h_0, N$ and $M$ and $\la$. Let $s=s_0$, we obtain  the following theorem. 

\begin{theorem}\label{10.14.T1}
	Under Assumption \ref{Assumption (H)}, for each weak solution $u$ of \eqref{10.10.1} with respect to $\la$, we have 
	\begin{equation}\label{10.14.10}
		\begin{split}
			\int_{\{x\in\Om\colon 0<\rho<r_0\}}u^2\df x\leq C\int_{\{x\in\Om\colon r_0<\rho<4r_0\}}u^2\df x, 
		\end{split}
	\end{equation}
	where the positive constants $C$ depend  only $\al, \be, r_0, R_0, h_0, N$ and $M$ and $\la$.
\end{theorem}

\begin{lemma}[{\cite[Theorem 5.6, p.~190]{Rousseau}}]\label{09.15.L1}
	Let $\Om$ be a connected open set in $\R^d$ and let $\om$ and $U$ be two open subsets of $\Om$ with $U\Subset\Om$. Let $C_0>0$. There exist $C>0$ and $\de\in (0,1)$ such that 
	\begin{equation*}
		\|u\|_{H^1(U)}\leq C\|u\|_{H^1(\Om)}^{1-\de}\left(\|f\|_{L^2(\Om)}+\|u\|_{L^2(\om)}\right)^\de, 
	\end{equation*}
	for $u\in H^2(\Om)$ satisfying 
	\begin{equation*}
		Pu=f+g, \mbox{ with } f\in L^2(\Om), \quad |g(x)|\leq C_0\left(|u(x)|+|Du(x)|\right) \mbox{ a.e.~in } \Om. 
	\end{equation*}
\end{lemma}

\begin{lemma}[{\cite[Lemma 5.10, p.~194]{Rousseau}}]\label{09.15.L2}
	Let $y\in \pt\Om$ and let $V$ be an open neighborhood of $y$ in $\R^d$ such that $\pt\Om$ is smooth in $V$. Let $C_0>0$. There exist $W\s V$ a neighborhood of $y$ in $\R^d$, $\e>0, \de\in (0,1)$, and $C>0$, such that we have 
	\begin{equation*}
		\|u\|_{H^1(W\cap \Om)}\leq C\|u\|_{H^1(\Om)}^{1-\de}\left(\|f\|_{L^2(\Om)}+\|u\|_{H^1(V\cap \Om_\e)}\right)^\de, \mbox{ with } \Om_\e=\{x\in\Om\colon \dist(x,\pt\Om)>\e\}
	\end{equation*}
	for $u\in H^2(\Om)$ satisfying $u|_{V\cap \pt\Om}=0$ and 
	\begin{equation*}
		Pu=f+g, \mbox{ with } f\in L^2(\Om), \quad |g(x)|\leq C_0(|u(x)|+|Du(x)|) \mbox{ a.e.~in }\Om. 
	\end{equation*}
\end{lemma}

\begin{theorem}\label{09.15.T1}
	Let $u$ be the weak solution of \eqref{10.10.1} with respect to $\la>0$. Let $\om\Subset \Om$ be an open subset. Then there exists a positive constant $C$, depending only on $\om, \al, \be, r_0, R_0, h_0, N$ and $M$ and $\la$, such that 
	\begin{equation}\label{09.15.1}
		\int_\Om u^2\df x\leq C\int_\om u^2\df x. 
	\end{equation}
\end{theorem}

\begin{proof}
	We divide the proof into two steps.

	{\it Step 1}. Propagation in the nondegenerate region.
	
	Set $m_\om:=\min_{\ol\om}\rho>0$, $a_\om:=\f{1}{2}m_\om$, and
	$$
	\Omega_{*,\om}:=\Omega\cap\{\rho>a_\om\}.
	$$
	Since $
	\rho^\alpha\geq a_\om^\alpha>0$ in $\Omega_{*,\om}$ and $\om\s\Omega_{*,\om}$, 
	the operator $
	P=-\Div(A\nabla\bullet)$ 
	is uniformly elliptic in $\Omega_{*,\om}$. Moreover, $u$ satisfies
	$$
	Pu=\la u \mbox{ in }\Omega_{*,\om}.
	$$ 
	Thus, Lemmas \ref{09.15.L1} and \ref{09.15.L2} can be applied in this nondegenerate region. By a finite chain of overlapping interior neighborhoods joining $\om$ to $\Omega\cap\{r_0<\rho<4r_0\}$, we obtain
	\begin{equation}\label{09.15.3}
		\|u\|_{L^2(\Omega\cap\{r_0<\rho<4r_0\})}
		\leq
		C\|u\|_{L^2(\om)},
	\end{equation}
	where $C$ depends on $\om, \al, \be, r_0, R_0, h_0, N$ and $M$ and $\la$
	
	{\it Step 2}. Propagation across the degenerate layer.
	
	By Theorem \ref{10.14.T1},
	\begin{equation}\label{09.15.4}
		\int_{\Om\cap\{0<\rho<r_0\}}u^2\,\df x
		\leq
		C
		\int_{\Om\cap\{r_0<\rho<4r_0\}}u^2\,\df x.
	\end{equation}
	Combining \eqref{09.15.3} and \eqref{09.15.4}, we obtain
	$$
	\int_{\Om}u^2\,\df x
	\leq
	C\int_{\om}u^2\,\df x.
	$$
	This proves \eqref{09.15.1}. We complete the proof of this theorem. 
\end{proof}

\begin{remark}
	For qualitative unique continuation of the degenerate elliptic equation
	\eqref{10.10.1}, Assumptions \ref{Assumption (H)} (1) and (3) are not
	required. Indeed, the proof relies only on Lemma \ref{10.10.L1},
	Remark \ref{10.10.R2}, and Lemma \ref{10.10.L2}, which establish the
	existence of a weak solution to \eqref{10.10.1}. Once the existence
	of a weak solution is established, the qualitative unique continuation
	property follows directly from the preceding arguments.
\end{remark}

\end{document}